\documentclass[letterpaper, paper,10pt]{AAS}		

\usepackage{bm}
\usepackage{amsmath}
\usepackage{amssymb}
\usepackage[colorlinks=true, pdfstartview=FitV, linkcolor=black, citecolor= black, urlcolor= black]{hyperref}
\usepackage{overcite}
\usepackage{footnpag} 
\usepackage{wrapfig}
\usepackage{graphicx}
\usepackage{subcaption}
\usepackage{float}
\usepackage{algorithm}
\usepackage{algpseudocode}
\usepackage{enumitem}
\usepackage{booktabs}
\usepackage{makecell}

\PaperNumber{25-831}

\begin{document}

\title{Long-Duration Station-Keeping Strategy for Cislunar Spacecraft Formations}

\author{Ethan Foss \thanks{PhD Student, Aeronautics and Astronautics, Stanford University.}, 
Yuji Takubo \thanks{PhD Student, Aeronautics and Astronautics, Stanford University.},
\ and Simone D'Amico \thanks{Associate Professor, Aeronautics and Astronautics, Stanford University.}
}

\maketitle{}

\begin{abstract}
This paper demonstrates a novel guidance and control strategy for cislunar near-rectilinear halo orbit formation-keeping applied to high-fidelity dynamics. Bounded relative motion is constructed about long-duration ephemeris trajectories with osculating invariant circles to form quasi-periodic relative orbits. State-of-the-art absolute control strategies are paired with a simple and effective relative control feedback law. Finally, a control barrier function is implemented to ensure recursively passively-safe bounded relative motion under feedback in the presence of possible missed maneuver events for the duration of the formation flight. The strategy is verified in high-fidelity simulation environments through Monte Carlo trials.
\end{abstract}

\section{Introduction}

Deployment of Distributed Space Systems (DSS) in cislunar orbits has been gaining interest following renewed attention in the exploration and economic activity of cislunar space.
In anticipation of the deployment of NASA's Gateway Station in a near-rectilinear halo orbit (NRHO) around the Earth-Moon $L_2$ libration point \cite{davisgateway}, development of DSS technologies, including swarming, formation flight, rendezvous, proximity operations, and docking (RPOD) requires further attention. 
In particular, long-term spacecraft swarming that performs bounded relative motion has the possibility to expand future space mission concepts.
However, due to the chaotic nature of the dynamical system around the Moon, traditional guidance, navigation, and control (GN\&C) algorithms that are tailored for the relative motion in perturbed Keplerian dynamics (e.g., Low Earth Orbit (LEO), Geostationary Orbit (GEO)) are not directly applicable to GN\&C systems for cislunar DSS. 
For example, the 9:2 Southern $L_2$ NRHO, which has received special interest for cislunar missions due to its increased stability, near-constant line-of-sight to Earth, and advantageous position for lunar and scientific operations, still offers numerous unexplored challenges for cislunar DSS, notably fully unstable relative motion.

In LEO, where the majority of distributed space systems have been deployed, the gravitational attraction of the Earth dominates other perturbative accelerations, such as those induced by gravitational harmonics, solar radiation pressure, third body gravitation, and drag. This favorable property means that trajectories in these regimes are not significantly perturbed, allowing for the applicability of low-fidelity models (i.e., the restricted two-body problem), for which semi-analytical solutions exist, and which offer numerous dynamical insights. Indeed, the existence of analytical relative motion models in near-earth orbits relies on the fact that absolute motion can be semi-analytically described \cite{koenig2017new}.
In contrast, low-fidelity models of spacecraft dynamics in cislunar space, such as the CR3BP, have not thus far produced similar levels of dynamical insight for absolute and relative motion. These models, which are non-linear and non-integrable, have not been found to admit analytical solutions. As a result, no exact and analytically available form for relative motion about periodic solutions has been found in these simplified models. Furthermore, simplified models perform significantly worse in describing true spacecraft motion in cislunar space, as a result of the stronger perturbative effects relative to the gravitational attraction of the primaries. For example, the assumption of circularity of the Earth-Moon orbit is a source of significant model error in the CR3BP, since the Moon's orbit of the Earth is non-negligibly eccentric and wobbling (e 	$\sim0.0255$ to $0.0775$) \cite{nasa_moonfactsheet}. Another large source of model error arises from solar gravitation. Inclusions of these effects, namely in the elliptical-restricted three-body problem (ER3BP) or the bi-circular restricted four-body problem (BCR4BP) \cite{parkassessment}, are still too oversimplified to describe orbital motion in cislunar space with sufficient levels of accuracy necessary for many applications.

Additionally, many orbits of interest in cislunar space are fully unstable (their state transition matrices contain eigenvalues outside of the unit circle) \cite{zimovan2020near}, in contrast to Keplerian orbits, which are only marginally unstable (their state transition matrices contain repeated eigenvalues on the unit circle) \cite{ALFRIEND201083}. This means that spacecraft will diverge from a reference trajectory at an exponential rate (in Keplerian orbits, trajectories diverge at a linear rate). As a result, any spacecraft in these orbits will need to frequently apply absolute and relative station-keeping maneuvers to maintain its orbit. In contrast, many spacecraft missions in LEO which do not demand high levels of absolute positioning accuracy can orbit the Earth for long-durations without requiring any absolute orbit station-keeping. 

As a result of the aforementioned perturbative nature of cislunar space, the design of absolute motion in cislunar orbits typically leverages procedures whereby trajectories are first designed in a low-fidelity model, which is easier to work with, provides dynamical insights, and admits periodic solutions, and then converged in a high-fidelity model, which includes perturbative effects and uses an ephemeris to accurately determine the locations of celestial bodies, to remove dynamical infeasibilities. This procedure is commonly used to generate long-duration baseline trajectories in the vicinity of periodic solutions of the CR3BP \cite{parkassessment,zimovanorbitgeneration}, and is similar to procedures used for some absolute station-keeping techniques in near-Earth orbit, such as for Tandem-X and Terra-SARx \cite{montenbruck2008navigation,krieger2007tandemx}. For many orbits of interest in the cislunar regime, this procedure serves as the basis for absolute orbit station-keeping. 

Several solution techniques have been proposed for performing absolute orbit control in cislunar orbits \cite{shirobokov2017survey, williams2023comparison}. Absolute orbit control is typically divided into long and short horizons \cite{davis2022orbit}, wherein long-horizon uncontrolled trajectories are generated in a high-fidelity ephemeris model (HFEM), using the aforementioned procedure of transitioning trajectories, and short-horizon correction maneuvers are performed to clean up dispersions along the long-duration reference trajectory. One of the most popular approaches, known as x-axis crossing \cite{davis2022orbit}, has been operated on the CAPSTONE mission \cite{cheetham2022capstone} and has been studied for use in Gateway station-keeping \cite{davis2022orbit}. It involves solving for the impulsive maneuver that targets an x-axis crossing of the reference trajectory several revolutions from the current time.
Other approaches, which all involve an HFEM reference, involve full-state feedback model predictive control \cite{shimane2025output}, chance-constrained feedback control design with missed-thrust events \cite{shrivastavamarkov}, and local eigenmotion control \cite{elangoeigenmotion}.

Relative motion control about periodic orbits of restricted multi-body problems has also been investigated and is typically described through linear dynamical models. In particular, Ref.~\citenum{elliott2022phd} derives a relative motion model that utilizes local toroidal coordinates to describe spacecraft relative motion and demonstrates its use for station-keeping under high-fidelity dynamics for up to 10 revolutions in the vicinity of a Sun-Earth $L_1$ libration point Halo orbit. Ref.~\citenum{takubo2025passively} extends these techniques by describing relative motion in a velocity-normal-binormal (VNB) frame and applying relevant optimal control techniques to local toroidal coordinates. The VNB frame is also utilized in Ref.~\citenum{vela2025modal}, in which approach trajectories are generated in modal coordinates by perturbing stable modes of the eigensystem. Other studies simply leverage cartesian coordinates, such as Ref.~\citenum{khoury_relative_2024}, which demonstrates continuous-thrust linear quadratic regulation in relative dynamics under relative navigation. Refs.~\citenum{takubo2025passively} and \citenum{takubo2025safeoptimalnspacecraftswarm} represent one of the few deployments of relative control of swarms under HFEM dynamics, though a key finding of these works is that the controller dynamics model error induces large dispersions and raises $\Delta V$ expenditure. 
 
Though absolute-relative formation-keeping in Keplerian regimes has been demonstrated, such as in the case of TANDEM-X \cite{krieger2007tandemx,montenbruck2008navigation}, no works have comprehensively formalized a strategy of formation-keeping in multi-body systems that considers the realistic conditions that present considerable difficulty to the success of formation flight in cislunar space. The HELIOSWARM mission \cite{helioswarm} involves 9 agent formation flight in a cislunar 2:1 mean motion resonance orbit. However, this orbit offers favorable stability properties and still leverages two-body relative motion models (i.e., the Yamanaka Ankersen STM \cite{yamanaka_new_2002}), meaning many of the station-keeping strategies leveraged by this mission are not extendible to a broad range of cislunar orbits. To the knowledge of the authors, relative orbit generation about NRHOs, which has been performed in simplified dynamical models like the CR3BP, ER3BP, and BCR4BP \cite{takubo2025safeoptimalnspacecraftswarm} has not been demonstrated about non-periodic high-fidelity reference trajectories, which allow for significantly lower $\Delta V$ expenditure over long durations. Additionally, perturbative effects to the relative motion induced by maneuver execution error and missed maneuver events in absolute control has not been considered in relative control contexts. Finally, many studies of spacecraft relative motion in cislunar space do not provide guarantees of passive safety under the possibility of missed maneuvers.

This work presents a novel strategy for cislunar formation flight that produces bounded absolute motion about near-rectilinear halo orbits and bounded relative motion of multi-agent formations with guarantees of passive safety. The strategy is simulated under realistic conditions with high-fidelity dynamics, maneuver execution error, missed maneuver events, and navigation uncertainty. By enforcing passive safety with a control barrier function (CBF), the strategy allows for decoupled and modularized computation of the absolute and relative control while efficiently and robustly enforcing formation safety. The resulting strategy produces low total $\Delta V$ expenditure on the order of 1 meter per second per agent per year, comparable to $\Delta V$ expenditures of single spacecraft station-keeping.

\section{Dynamical Models}

Absolute motion of spacecraft in the cislunar regime can be described at varying levels of fidelity. This work leverages the CR3BP and a high-fidelity ephemeris model (HFEM) to describe spacecraft motion in this regime. Moreover, relative motion models can be constructed about reference trajectories of the two models.

\subsection{CR3BP Dynamics}

The CR3BP is frequently used to model spacecraft motion in cislunar space. It makes several assumptions to simplify resulting equations of motion, namely circular motion of the primaries and negligible third body (spacecraft) mass. Additionally, states are normalized by characteristic lengths and times and expressed in a synodic frame. The resulting equation of motion is
\begin{equation}
    \ddot{\bm{r}} = -2 \hat{\bm{z}} \times \bm{r} + \nabla \Omega.
\end{equation}
\noindent $\Omega$ is the pseudo-potential of the CR3BP, given by
\begin{equation}
    \Omega = \frac{1-\mu}{r_1} + \frac{\mu}{r_2} + \frac{x^2 + y^2}{2}
\end{equation}
\noindent where $\mu$ is a mass ratio constant associated with the two primaries and $r_1$ and $r_2$ denote distances to the first and second primary respectively.

\subsection{HFEM Dynamics}

High-fidelity dynamical models are used to construct long-duration trajectories and verify resulting control schemes. The equations of motion are expressed in a Moon-fixed normalized rotating frame \cite{parkassessment}. This results in the equation of motion
\begin{equation}
    \ddot{\bm{r}} = M^{\top} \left( -\ddot{M} \bm{r} - 2 \dot{M} \dot{\bm{r}} + \ddot{\bm{R}} - \ddot{\bm{B}} \right).
\end{equation}
\noindent Here, $M$ is a direction cosine matrix constructed from ephemeris positions and velocities of the two primaries that describes the rotation from the inertial frame to the synodic frame. $\ddot{\bm{B}}$ is the acceleration of the Moon in the inertial frame, and $\ddot{\bm{R}}$ is the acceleration of the spacecraft in the inertial frame, which can be modeled with gravitational accelerations from the primaries (e.g. Earth, Moon), additional gravitational perturbations (e.g. Sun, Venus), and possibly solar radiation pressure and gravitational harmonics. Relevant ephemeris data is obtained from SPICE \cite{actonspice}.

\subsection{Relative Motion Dynamics}

Linear relative dynamics can be produced by taking a first order approximation of the absolute dynamics about a reference trajectory. Let $\bm{\phi}(t)$ be a reference trajectory which is a solution to the set of equations of motion $\dot{\bm{x}}(t) = \bm{f}(\bm{x}(t),t)$. If $\delta \bm{x}(t) = [ \delta \bm{r}(t)^\intercal, \delta \bm{v}(t)^\intercal]^\intercal$ is defined with respect to $\bm{\phi}(t)$ then 
\begin{equation}
    \delta \dot{\bm{x}}(t) = \bm{f}(\bm{\phi}(t) + \delta \bm{x}(t)) - \bm{f}(\bm{\phi}(t)).
\end{equation}
\noindent Taking the first-order variation about $\delta \bm{x}(t) = 0$ gives
\begin{equation}
    \delta \dot{\bm{x}}(t) \approx A(t) \delta \bm{x}(t)
\end{equation}
\noindent where $A(t) \triangleq \frac{\partial \bm{f}}{\partial \bm{x}} \Big|_{\bm{\phi}(t)}$ is the Jacobian of the absolute dynamics evaluated at the reference trajectory. The state transition matrix (STM) can also be constructed to provide a discrete form of the relative motion model with 
\begin{equation}
    \dot{\Phi}(t,t_0) = A(t) \Phi(t,t_0), \quad \Phi(t_0,t_0) = I.
\end{equation}

This approach to the construction of the relative motion dynamics is agnostic to the absolute dynamics, meaning this approach is valid for both the CR3BP and HFEM dynamics. However, this relative motion must be numerically integrated, and has no analytical solution, unlike for Keplerian regimes \cite{ALFRIEND201083}. 

\subsection{Coordinate Frames} 

Expressing the relative motion in different coordinate frames can often provide additional insight into underlying dynamical structure. Relative motion in Keplerian regimes is often expressed in LVLH frames. A popular approach for cislunar regimes is to express relative motion in a velocity-normal-binormal (VNB) frame (see Refs.~\citenum{takubo2025passively,vela2025modal}), where axes are aligned with the along track velocity and moon angular momentum directions (also sometimes referred to as a velocity-normal-conormal (VNC) frame \cite{helioswarm}). Let the reference trajectory be decomposed into its positions and velocities, i.e. $\bm{\phi}(t) = [\bm{r}^\intercal_{syn}(t), \bm{v}^\intercal_{syn}(t)]^\intercal$ and let $\bm{r}_{M,syn}$ be the position of the moon with respect to the system origin. A transformation from a relative position in the synodic frame to the VNB frame can be performed with 
\begin{equation}
\begin{gathered}
    \delta \bm{r}_{syn}(t) = R(t) \delta \bm{r}_{VNB}(t), \quad R(t) = \begin{bmatrix}
        \hat{\bm{x}} & \hat{\bm{y}} & \hat{\bm{z}}
    \end{bmatrix}, \\ \quad \hat{\bm{x}} \parallel \bm{v}_{syn}(t), \quad \hat{\bm{y}} \parallel \left( \bm{v}_{syn}(t) \times (\bm{r}_{syn}(t) - \bm{r}_{M,syn}) \right), \quad \hat{\bm{z}} = \hat{\bm{x}} \times \hat{\bm{y}}.
    \end{gathered}
\end{equation}




\section{Absolute and Relative Motion Generation}

\subsection{Absolute Trajectory}

\begin{wrapfigure}{r}{0.25\textwidth}  
\vspace{-2cm}
  \centering
  \includegraphics[width=0.23\textwidth]{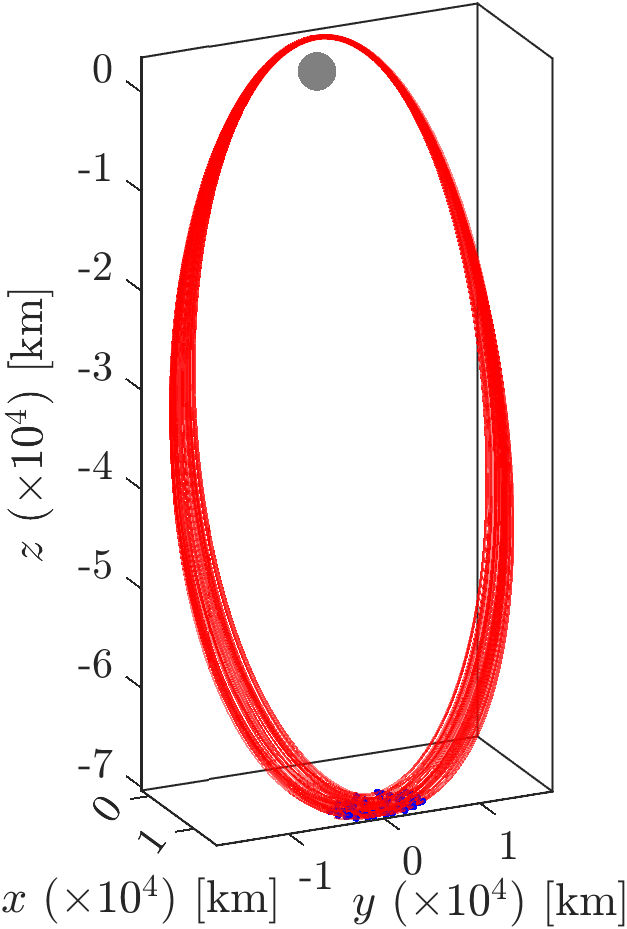}
  \caption{\small HFEM trajectory near the 9:2 NRHO, generated at an epoch of 2024-10-29 12:00:00 TDB and with Sun gravitational perturbations. Patch points (blue) are at apolune.}
  \label{fig:Ephem25to6}
  \vspace{-1.5cm}
\end{wrapfigure}

Performing long-duration formation-keeping relies upon the generation of reference trajectories to be tracked by a spacecraft. For absolute motion, this is typically done by taking a periodic orbit in a lower fidelity model (e.g., CR3BP), and performing numerical continuation to a HFEM \cite{davis2022orbit,parkassessment}. Ephemeris trajectories are generated by stacking the apoapsis node of a CR3BP orbit for the desired number of revolutions and converging the initial guess in the HFEM with a multiple shooting scheme. Fig.~\ref{fig:Ephem25to6} shows an example of a 50 revolution ephemeris trajectory generated in the vicinity of a 9:2 NRHO from the CR3BP. For more details on trajectory generation in ephemeris models, see Refs.~\citenum{davis2022orbit,zimovanorbitgeneration}.

\subsection{Relative Orbit Reference Generation}

\begin{figure}[h]
    \centering
    \includegraphics[width=0.6\linewidth]{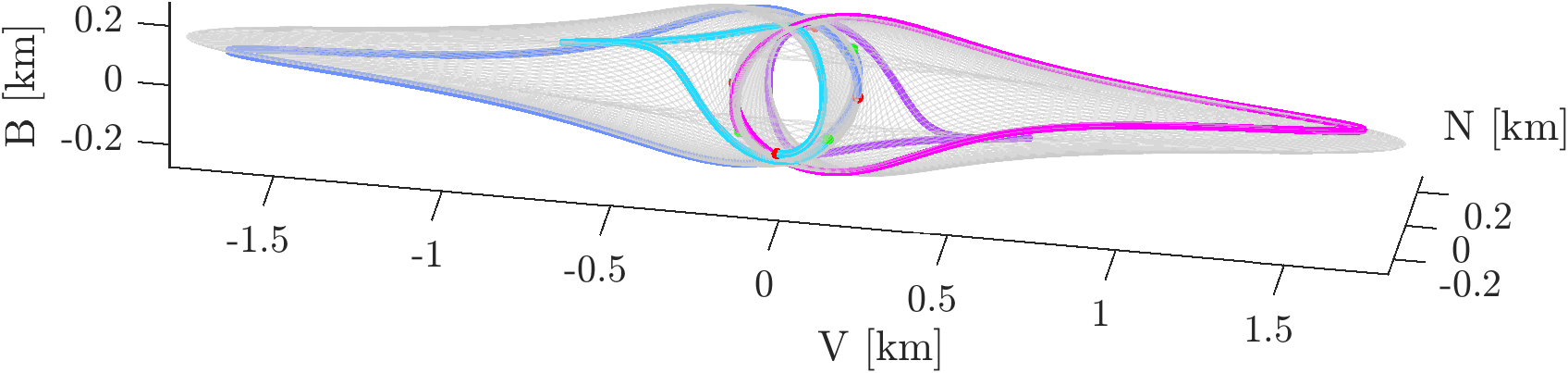}
    \caption{Relative motion in the CR3BP generated from a quasi-periodic torus about the 9:2 NRHO}
    \label{fig:CR3BPQPO}
\end{figure}

The generation of bounded relative motion between multiple spacecraft in Keplerian orbits simply relies on the satisfaction of the energy matching condition, that is, the orbits of interest have equal semi-major axes \cite{schaubspacecraftformationflying}. For cislunar orbits, generating bounded relative motion can be more challenging, due to the numerical nature of the trajectories. Since periodic orbits exist in one-parameter families in the CR3BP, neighboring periodic orbits will not produce bounded relative motion in general. However, quasi-periodic tori, which exist about periodic solutions, can produce bounded relative motion so long as its period is equivalent with that of the periodic orbit. Because of this, quasi-periodic tori have served as the foundation for generating bounded relative motion in cislunar space \cite{elliott2022describing, baresi2017spacecraft, takubo2025passively}. In close proximity to a periodic orbit, these quasi-periodic orbits exist as a linear combination of the eigenvectors associated with unimodular eigenvalues of the Monodromy matrix. The real and imaginary components of these eigenvalues span an invariant circle, which when propagated for a full revolution, repeats on itself, forming a quasi-periodic relative orbit (QPRO). A QPRO about the 9:2 NRHO in the CR3BP along with a selection of relative trajectories on its surface are displayed in Fig.~\ref{fig:CR3BPQPO} in the VNB frame.

Unlike the CR3BP, ephemeris trajectories are non-periodic, meaning globally bounded relative motion does not exist in general. However, despite the non-periodicity of ephemeris trajectories, their associated state transition matrices do exhibit central modes that are close to unimodular after integer revolutions, which can be perturbed to produce locally bounded relative motion about the trajectory for a finite duration, similar to the procedure used to generate quasi-periodic tori \cite{olikara2016computation,foss2024thesis}. If $\bm{\phi}(t)$ represents an ephemeris trajectory $\forall t \in T$, and $\Phi(t,t_0)$ is the associated state transition matrix of the trajectory, then a point in time $t_1 \in T$ can be identified for which a nearly central mode exists. The eigenvalue and eigenvector associated with this central mode can be denoted by $e_r\pm e_i i$ and $\bm{w}_r \pm \bm{w}_i i$. In dynamics that exhibit periodicity (e.g., CR3BP), this central mode is unimodular for the Monodromy matrix and produces an invariant circle. In ephemeris dynamics, an osculating invariant circle can be constructed in a similar manner with
\begin{equation}
    \bm{u}(t_0,\theta_1) = K ( \cos(\theta_1) \bm{w}_r -\sin(\theta_1) \bm{w}_i )
\end{equation}
\noindent where $K$ is the radius of the invariant circle and $\theta_{1}$ is the initial phase angle. After the full duration of relative motion, the final invariant circle is
\begin{equation}
\begin{aligned}
    \bm{u}(t_1,\theta_1) &= \Phi(t_1,t_0) \bm{u}(t_0,\theta_1) \\ &= K \left[(e_r \cos(\theta_1) - e_i \sin(\theta_1)) \bm{w}_r - (e_r \sin(\theta_1) + e_i \cos(\theta_1))\bm{w}_i \right]
    \end{aligned}
\end{equation}

\noindent It is worth noting that bounded relative motion generated from these osculating invariant circles is only guaranteed for $t \in [t_0,t_1]$. Moreover, the extent of the deviation in this time frame cannot necessarily be bounded below some finite value. Finally, the relative motion is only studied under linear dynamics, but for fully nonlinear relative motion, it is expected that large deviations from the nominal trajectory will lead to dispersions.

\begin{figure}[H]
    \centering

    \begin{subfigure}[t]{0.49\textwidth}
        \centering
        \includegraphics[width=\textwidth]{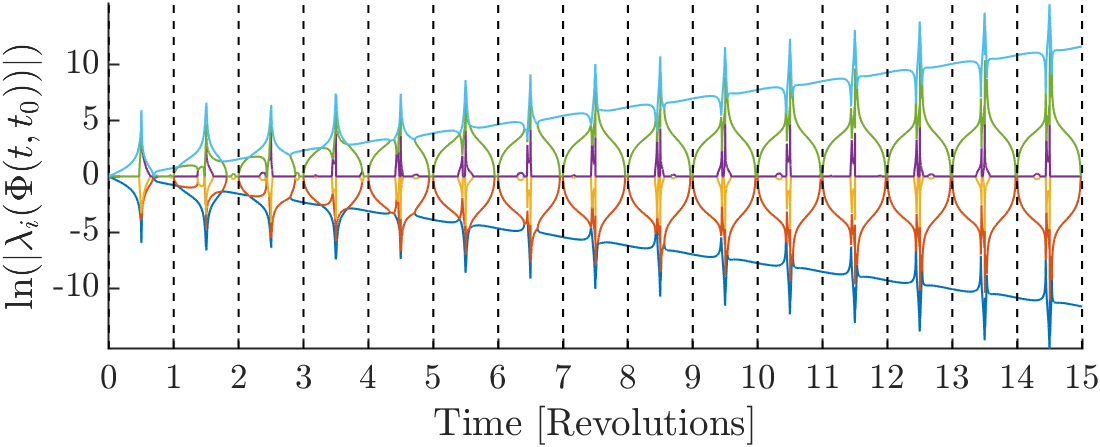}
        \caption{CR3BP}
    \end{subfigure}
    \hfill
    \begin{subfigure}[t]{0.49\textwidth}
        \centering
        \includegraphics[width=\textwidth]{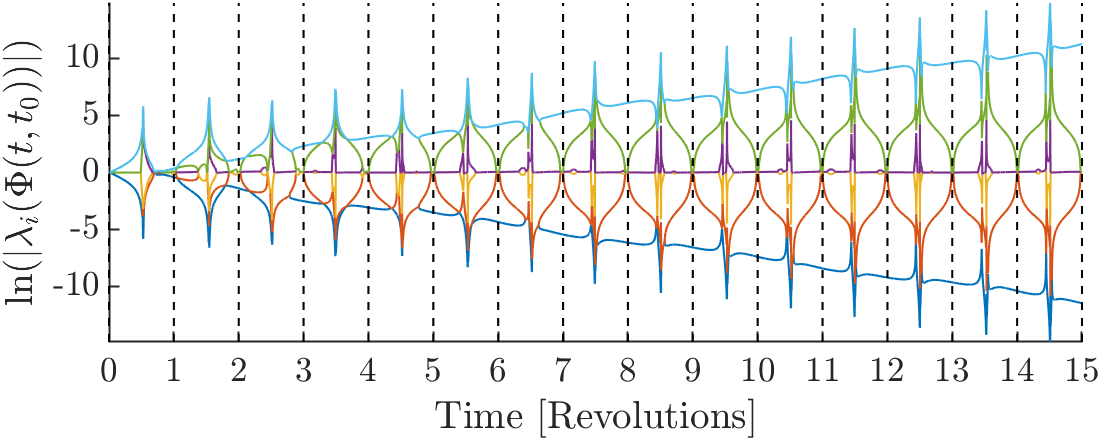}
        \caption{HFEM}
    \end{subfigure}
    \caption{Growth of logarithm of eigenvalues of STMs over time}
    \label{fig:EigenGrowth}
\end{figure}

To compare the eigenstructure between the CR3BP and HFEM, from which QPROs are generated, Fig.~\ref{fig:EigenGrowth} visualizes the growth over time of the eigenvalues associated with the 9:2 NRHO of the CR3BP and the ephemeris trajectory in Fig.~\ref{fig:Ephem25to6}. Evidently, the eigenstructure of the HFEM trajectory closely resembles the eigenstructure of the CR3BP, which is repeating, allowing for the construction of QPROs about HFEM trajectories over long durations. Unfortunately, there are still limitations to this approach of generating bounded relative motion. First, the duration of an osculating QPRO is limited by the fact that over sufficiently long durations, unimodular eigenvalues cease to exist and the condition number of the STM becomes unmanageably large. Second, even if unimodular eigenvalues exist, the QPRO generated from them may exhibit very large range ratios, that is, the maximum deviation from the origin may be significantly larger than the minimum deviation from the origin, which is generally undesirable for formation flight. This means that in contrast to CR3BP relative motion, for which QPROs exist for infinite duration, HFEM relative motion must be periodically regenerated. 

\begin{figure}[htbp]
  \centering

  \begin{subfigure}{\textwidth}
    \centering
    \begin{subfigure}{0.6\textwidth}
      \includegraphics[width=\linewidth]{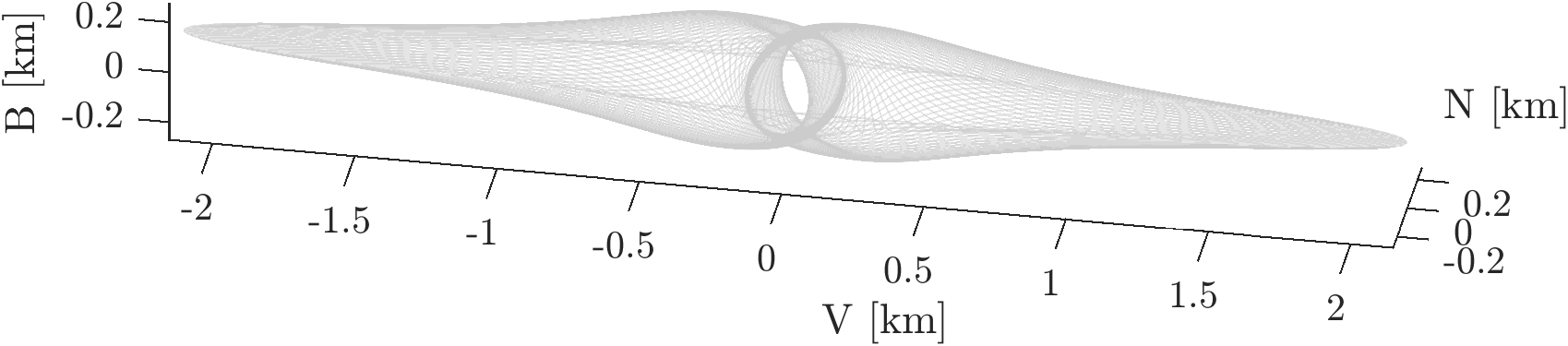}
    \end{subfigure}
    \hfill
    \begin{subfigure}{0.35\textwidth}
      \includegraphics[width=\linewidth]{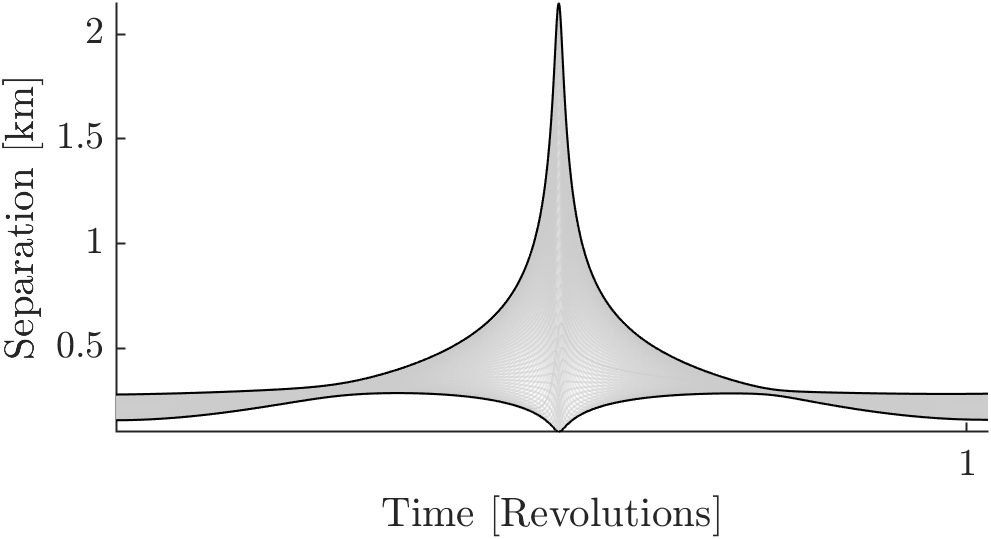}
    \end{subfigure}
    \caption{One Revolution}
  \end{subfigure}

  \vspace{0cm} 

  \begin{subfigure}{\textwidth}
    \centering
    \begin{subfigure}{0.6\textwidth}
      \includegraphics[width=\linewidth]{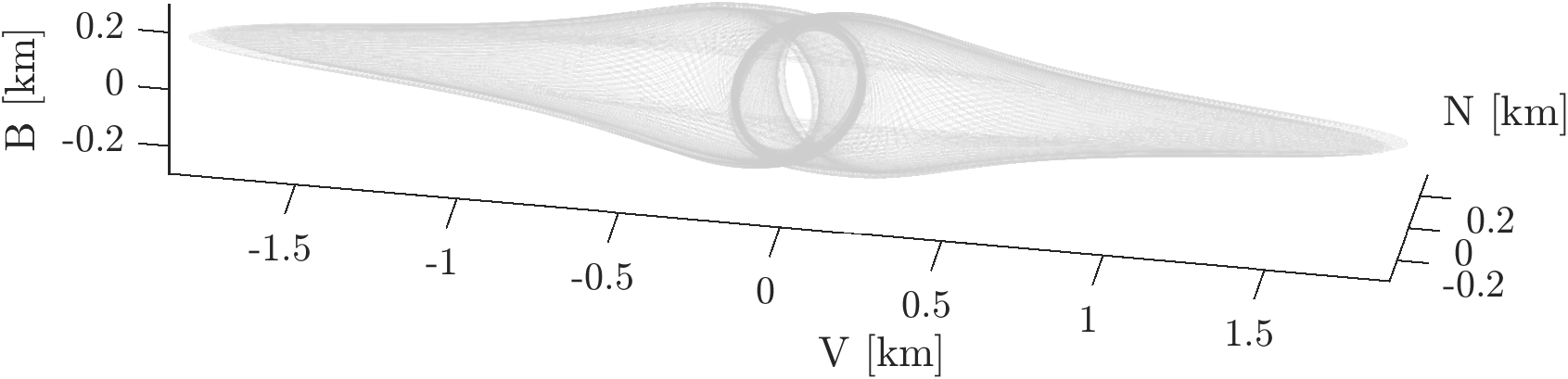}
    \end{subfigure}
    \hfill
    \begin{subfigure}{0.35\textwidth}
      \includegraphics[width=\linewidth]{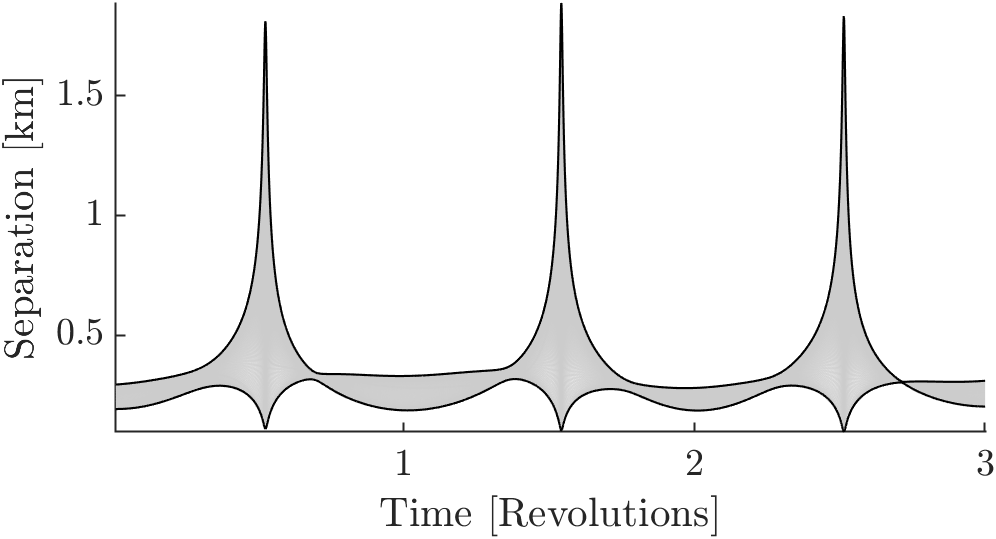}
    \end{subfigure}
    \caption{Three Revolutions}
  \end{subfigure}

    \vspace{0cm} 

  \begin{subfigure}{\textwidth}
    \centering
    \begin{subfigure}{0.6\textwidth}
      \includegraphics[width=\linewidth]{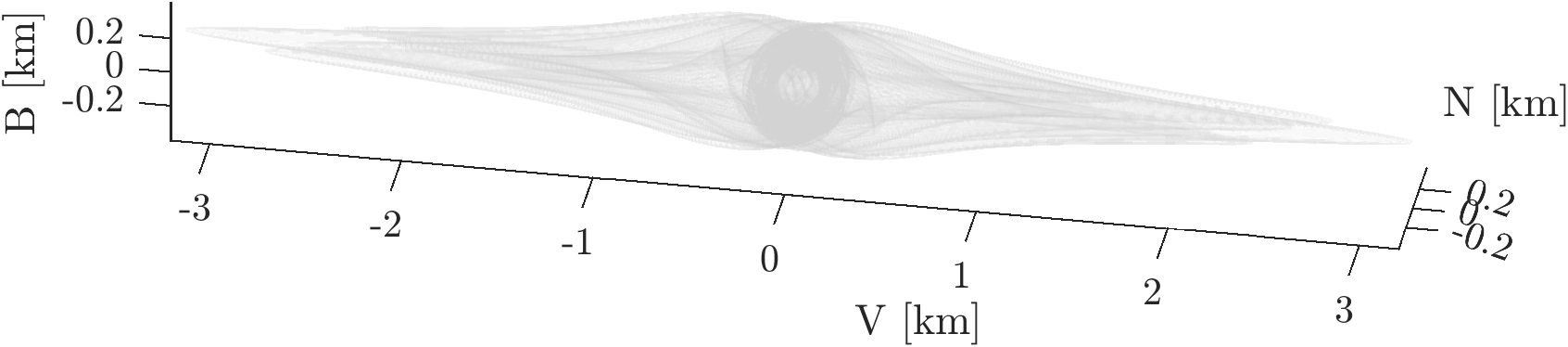}
    \end{subfigure}
    \hfill
    \begin{subfigure}{0.35\textwidth}
      \includegraphics[width=\linewidth]{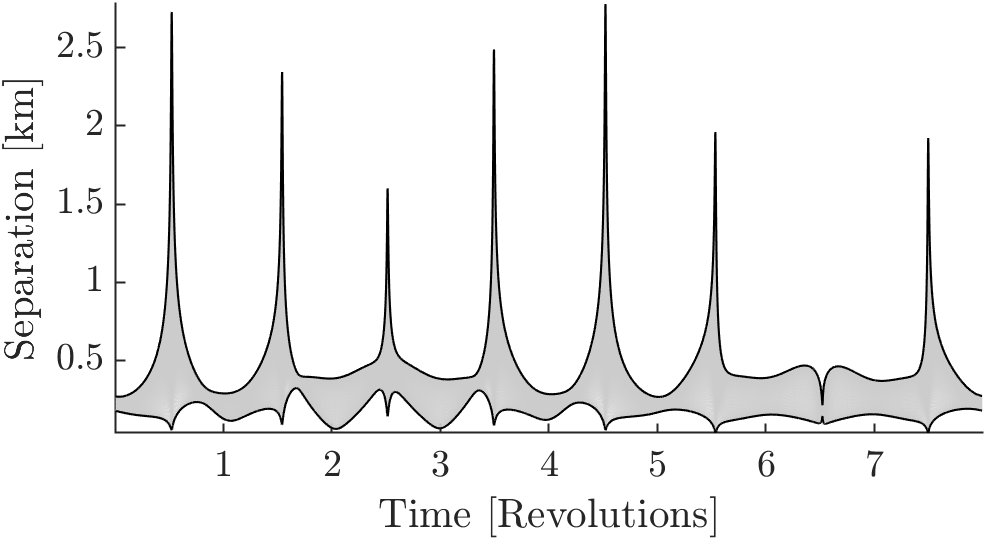}
    \end{subfigure}
    \caption{Seven Revolutions}
  \end{subfigure}

    \vspace{0cm} 

  \begin{subfigure}{\textwidth}
    \centering
    \begin{subfigure}{0.6\textwidth}
      \includegraphics[width=\linewidth]{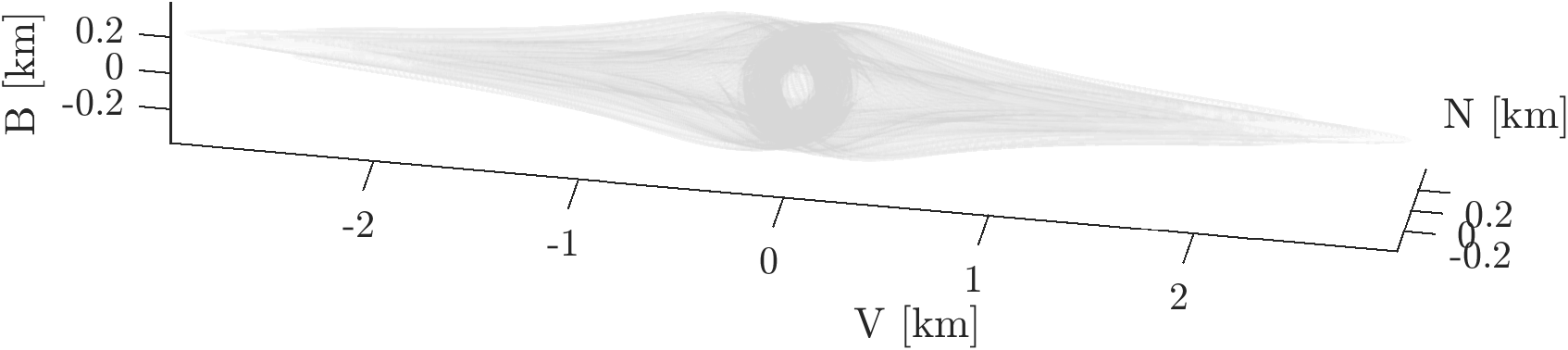}
    \end{subfigure}
    \hfill
    \begin{subfigure}{0.35\textwidth}
      \includegraphics[width=\linewidth]{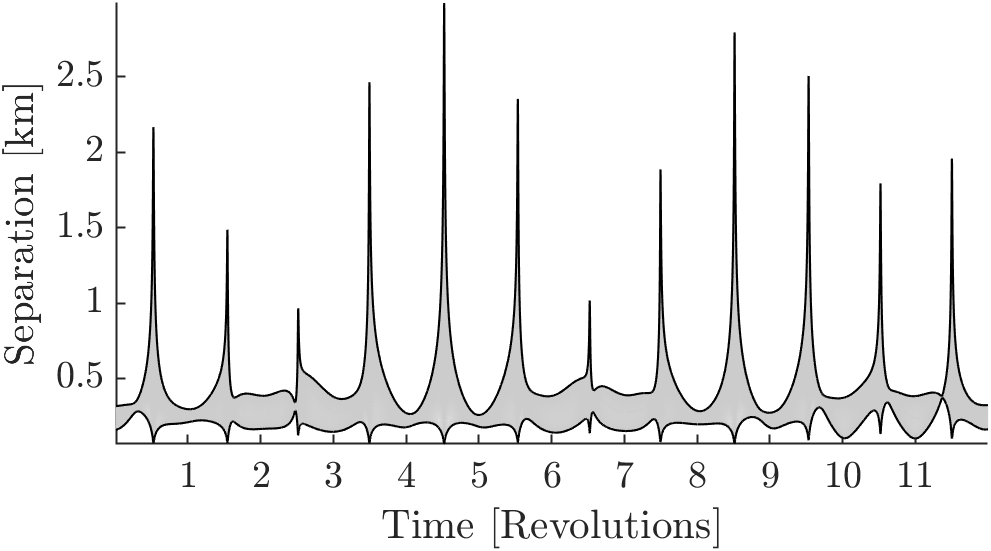}
    \end{subfigure}
    \caption{Twelve Revolutions}
  \end{subfigure}

    \vspace{0cm} 

  \begin{subfigure}{\textwidth}
    \centering
    \begin{subfigure}{0.6\textwidth}
      \includegraphics[width=\linewidth]{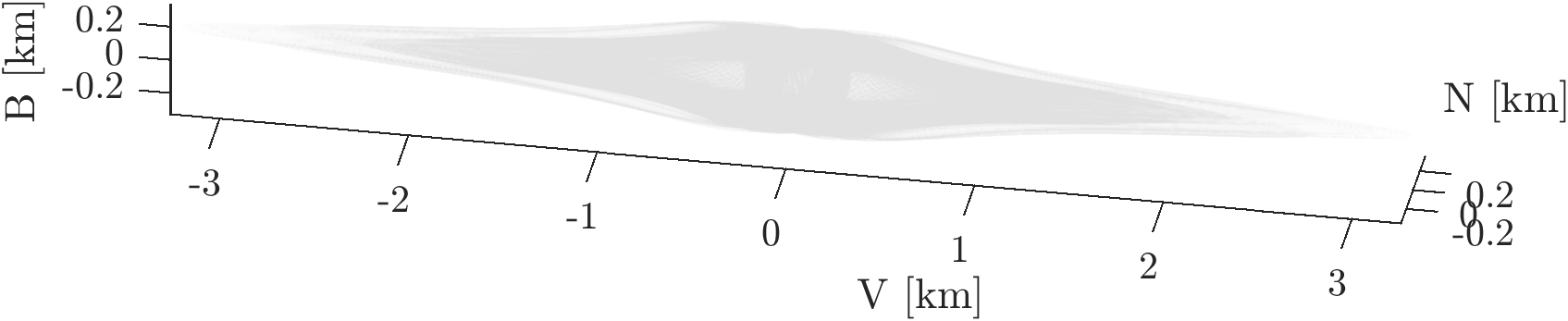}
    \end{subfigure}
    \hfill
    \begin{subfigure}{0.35\textwidth}
      \includegraphics[width=\linewidth]{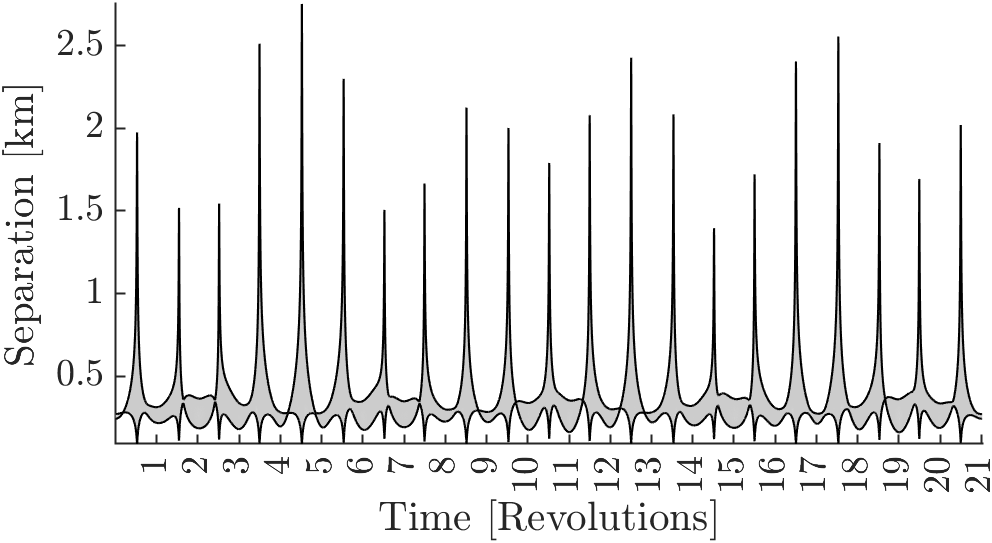}
    \end{subfigure}
    \caption{21 Revolutions}
  \end{subfigure}

  \caption{QPROs about a 9:2 HFEM trajectory at various revolution counts}
  \label{fig:EphemQPROs}
\end{figure}

The 9:2 NRHO of the CR3BP and its associated HFEM trajectories admit two sets of (quasi-)unimodular eigenvalues. One of these corresponds to the trivial along-track separated relative motion, and the other corresponds to the QPRO. In theory, a combination of these two modes can be perturbed to produce bounded relative motion, though this paper will focus on the modes associated with the QPROs. Fig.~\ref{fig:EphemQPROs} displays a selection of QPROs generated at varying integer revolution counts. The QPROs closely resemble the geometry of the CR3BP QPROs (Fig.~\ref{fig:CR3BPQPO}) but gradually lose their shape and produce larger separations with higher integer revolutions. There exist two degrees of freedom for spacecraft placement along QPROs, with all spacecraft being coplanar. An example of how a formation geometry can be constructed is visualized in Fig.~\ref{fig:InvariantCircleFormations}.  

Twelve revolutions of the 9:2 HFEM NRHO, or about 80 days of flight, corresponding to an STM condition number of $1e8$, balances the maintenance of QPRO geometry and frequency of reference switching. Unfortunately, there is no good technique for minimizing the range ratio of a QPRO. As such, employing a brute force approach, whereby a selection of time values are chosen for the final STM and tested for range ratio, is the procedure used for producing favorable QPRO geometry about HFEM trajectories. The existence of QPROs in high-fidelity dynamics is nonetheless an important result for cislunar formation flight, as it allows for the generation of long-duration bounded relative motion in realistic and relevant formation flight scenarios.

\subsection{Local Toroidal Coordinates About HFEM Trajectories}

\begin{figure}
    \centering
    \includegraphics[width=0.5\linewidth]{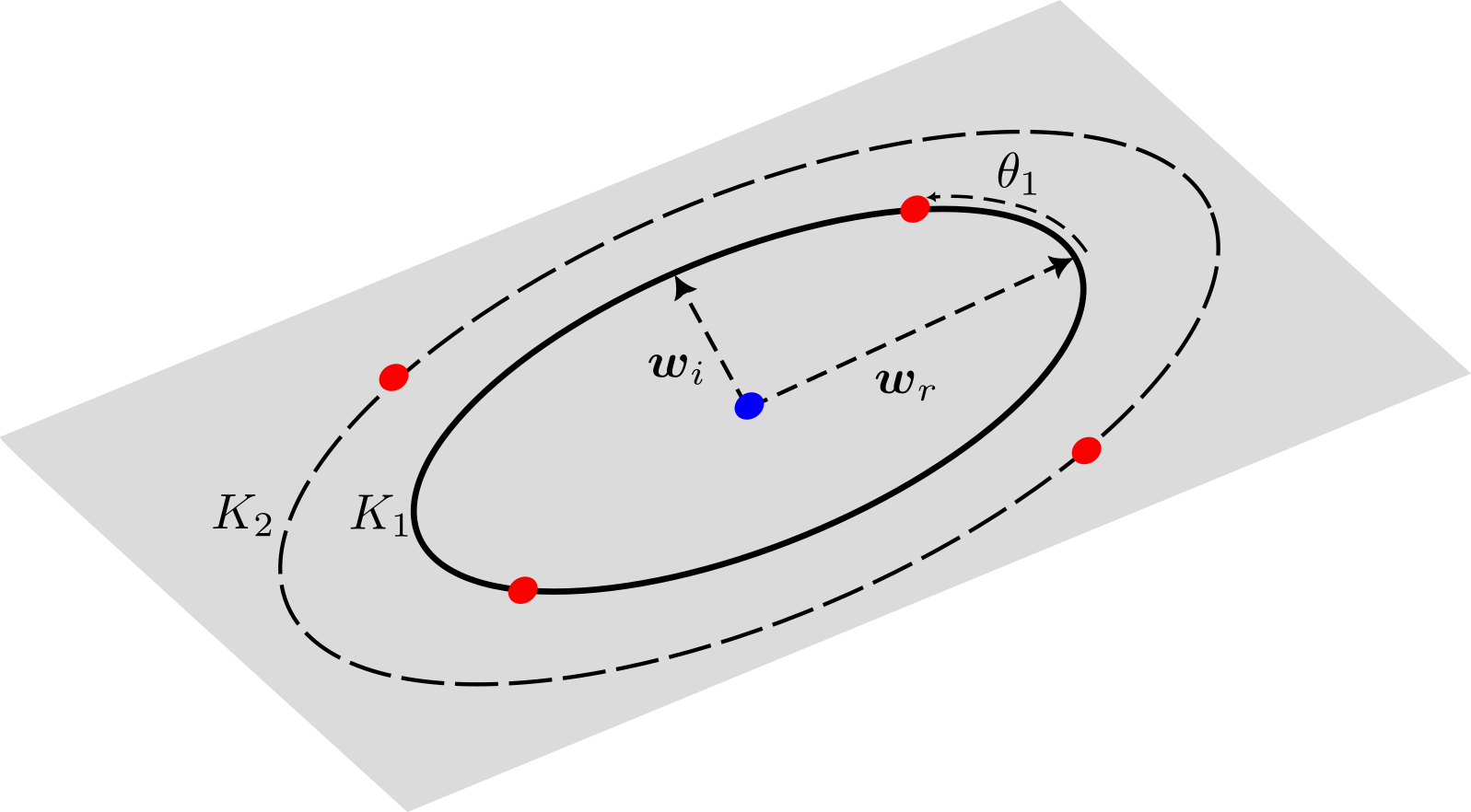}
    \hfill
    \includegraphics[width=.3\linewidth]{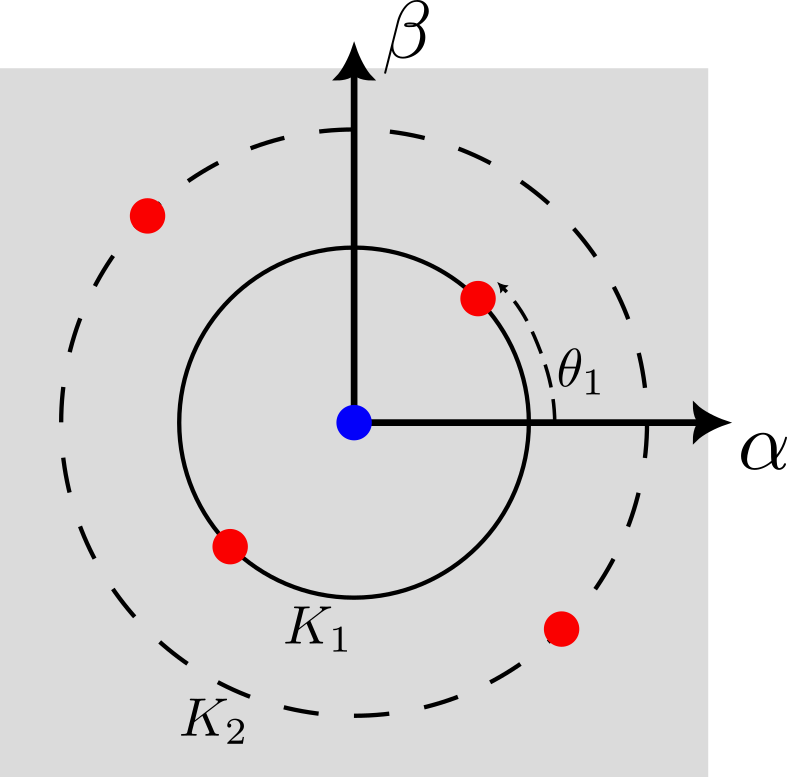}
    \caption{Formation configuration on an osculating invariant circle. Spacecraft can be separated by phase angle ($\theta_1$) or radius ($K$), and equivalently represented in local toroidal coordinates ($\alpha$,$\beta$). }
    \label{fig:InvariantCircleFormations}
\end{figure}

Local toroidal coordinates (LTCs) offer a convenient way to visualize, represent, and encode equilibrium in relative orbital geometries in restricted three-body problem orbits \cite{elliott2022describing,elliott2022phd}. As a result of the procedure used to identify finite-duration QPROs about HFEM trajectories, it is possible to leverage LTCs for the duration of the QPRO despite the fact that the underlying reference is not periodic. Let  $e_r+e_ii$ represent the quasi-unimodular eigenvalues and $\bm{w}_{r,0} +\bm{w}_{i,0} i$ represent their associated eigenvectors determined from a HFEM STM $\Phi(t_1,t_0)$. The mapping of these eigenvectors to any point in time along the reference trajectory is given by $\bm{w}_r(t) = \Phi(t,t_0) \bm{w}_{r,0}$ and  $\bm{w}_i(t) = \Phi(t,t_0) \bm{w}_{i,0}$. Let the eigenvectors be decomposed into position and velocity components with $\bm{w}_r = \begin{bmatrix} \bm{r}_r^\intercal , \bm{v}_r^\intercal \end{bmatrix}^\intercal$ and $\bm{w}_i = \begin{bmatrix} \bm{r}_i^\intercal , \bm{v}_i^\intercal \end{bmatrix}^\intercal$. Local toroidal coordinates are defined as $\bm{\zeta}\triangleq \begin{bmatrix} \bm{\zeta}_r^\intercal, \bm{\zeta}_v^\intercal\end{bmatrix}^\intercal$, $\bm{\zeta}_r \triangleq [\alpha, \beta, h]^\intercal$, $\bm{\zeta}_v \triangleq [\dot{\alpha}, \dot{\beta}, \dot{h}]^\intercal$ and can be related to a relative position and velocity with respect to the reference with
\begin{equation}
\begin{gathered}
    \delta \bm{r}(t) = T(t) \bm{\zeta}_r,  \quad T(t) = \begin{bmatrix}
        \bm{r}_r(t) & \bm{r}_i(t) & \hat{\bm{n}}(t)
    \end{bmatrix}, \quad \hat{\bm{n}}(t) \parallel \bm{r}_r(t) \times \bm{r}_i(t) \\
    \delta \bm{v}(t) = \dot{T}(t) \bm{\zeta}_r + T(t) \bm{\zeta}_v,  \quad \dot{T}(t) = \begin{bmatrix}
        \bm{v}_r(t) & \bm{v}_i(t) & \dot{\hat{\bm{n}}}(t)
    \end{bmatrix}.
\end{gathered}
\end{equation}
\noindent If the LTCs of a spacecraft are such that $h=0$ and $\bm{\zeta}_v = \bm{0}$, then the point lies on the surface of a locally invariant torus. Moreover, $\alpha$ and $\beta$ describe the location on the invariant circle, as is visualized in Fig.~\ref{fig:InvariantCircleFormations}.

\section{Absolute and Relative Station-Keeping Control}

Due to the instability of NRHOs in cislunar space, small perturbations and dispersions will lead to exponential divergence from reference trajectories. As such, station-keeping maneuvers must be performed periodically to keep a spacecraft near its reference. In comparison to Keplerian regimes, cislunar space poses several additional challenges to formation-keeping control. Namely, the presence of strong sources of perturbation, the possibility of missed maneuver events (MMEs), and the risk of collisions render the absolute and relative station-keeping highly coupled. This paper proposes a strategy that nonetheless allows for separate computation of absolute and relative station-keeping maneuvers by utilizing a control barrier function (CBF) safety filter to enforce keep out zones, providing guarantees of safe operation under uncertainty. Absolute control, in these contexts, refers to the station-keeping of the absolute state of the virtual chief of the formation to the absolute HFEM trajectory (e.g., from Fig.~\ref{fig:Ephem25to6}), whereas relative control refers to the station-keeping of the individual agents to their prescribed relative orbital geometry, as determined from a reference QPRO.

\begin{figure}
    \centering
    \includegraphics[width=0.8\linewidth]{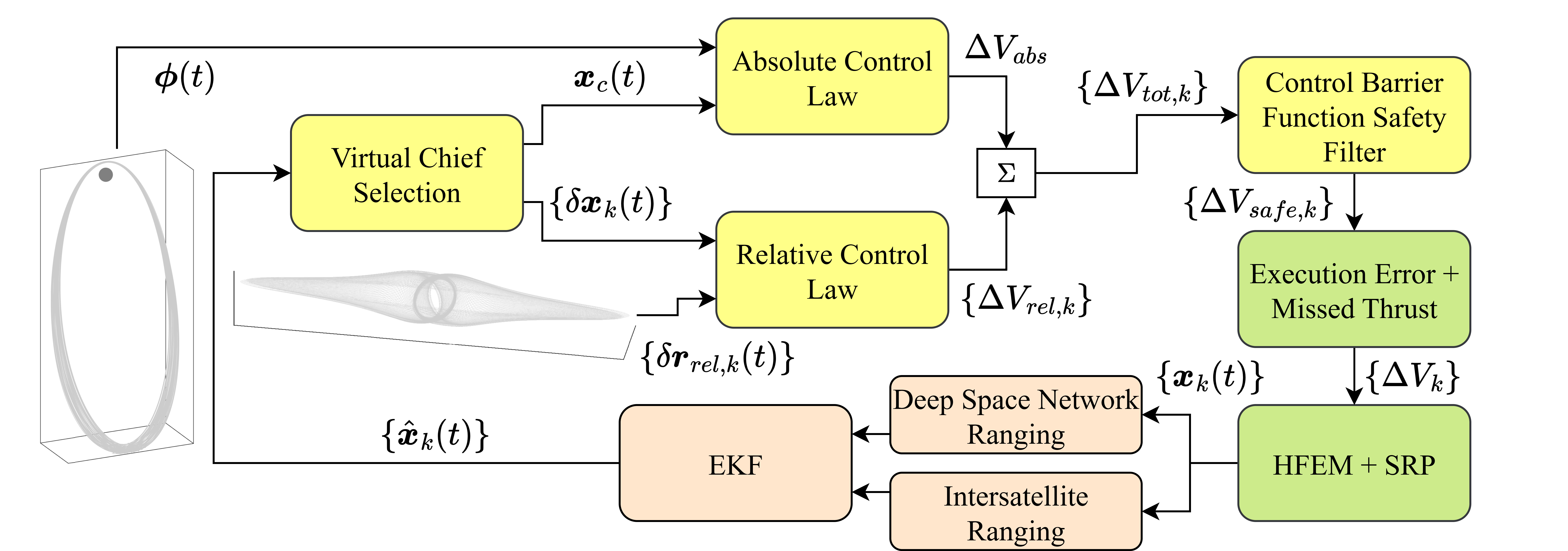}
    \caption{Proposed GN\&C architecture}
    \label{fig:GNCDiagram}
\end{figure}

For a formation consisting of $M$ agents, the proposed scheme computes an absolute station-keeping maneuver for a virtual chief that targets a perilune position at several revolutions downstream along an HFEM trajectory. This maneuver is then applied to all agents of the formation. Relative control is performed similarly, targeting a relative position with respect to the virtual chief along a pregenerated QPRO. The combination of absolute and relative maneuvers is then filtered by the CBF, which guarantees passive safety. This process is summarized in the GN\&C block diagram in Fig.~\ref{fig:GNCDiagram}.

\subsection{Absolute Control}


To compute absolute station-keeping maneuvers that keep the formation in the vicinity of a reference, this work makes use of an MPC scheme similar to the one proposed in Ref.~\citenum{shimane2025output}. However, the scheme proposed here instead computes just one maneuver (instead of two, as is done in Ref.~\citenum{shimane2025output}) that targets a position ellipsoid at perilune crossing six revolutions downstream along the reference. This is akin to the popularly employed x-axis crossing scheme \cite{davis2022orbit,Guzzetti2017}, which involves computing one maneuver to target the state of the reference trajectory at six revolutions downstream, which has been found to produce favorable $\Delta V$ expenditures. Given a reference trajectory $\bm{\phi}(t)$, an initial condition $\bm{x}_0$ and initial time $t_0$, the $\Delta V_{abs}$ is solved for with
\begin{equation}
\begin{aligned}
    \min_{\Delta V_{abs}} \quad & \| \Delta V_{abs}\|_2 \\
    \text{s.t.} \quad & \bm{x}_f = \bm{x}_0 + B \Delta V_{abs} + \int_{t_0}^{t_f} \bm{f}(\bm{x},t) \text{d} t \\
                      & \|[I \ \bm{0}] (\bm{x}_f-\bm{\phi}(t_f))\|_2 \le \epsilon,
\end{aligned}
\end{equation}
\noindent where $B =\begin{bmatrix} \bm{0} & I \end{bmatrix}^\intercal$, $t_f$ is the time corresponding to perilune crossing at six revolutions downstream from the current time, and $\bm{f}$ is a function representing the HFEM equations of motion of the system. Given the nonlinear nature of the optimization problem, it is solved using sequential convex programming (SCP) by enforcing a trust region about the solution variables, augmenting the dynamics constraint with a slack variable, and penalizing the $l_2$ norm of the slack in the objective. The only solution variables in the SCP problem are the initial condition dynamic slack variable $\bm{\nu} \in \mathbb{R}^6$, which is augmented to the objective with an $l_2$-norm penalty, and the $\Delta V_{abs} \in \mathbb{R}^3$, for which a trust region of 1 cm/s is enforced. The final position excursion ellipsoid is set to $\epsilon=$ 25 km. The optimization routine should be theoretically always feasible \cite{shimane2025output}; and an appropriately applied SCP solution routine robustly and reliably converges \cite{mao2019successiveconvexificationsuperlinearlyconvergent}. 

\begin{figure}[htbp]
  \centering
  \begin{minipage}{0.4\textwidth}
    \centering
    \includegraphics[width=\linewidth]{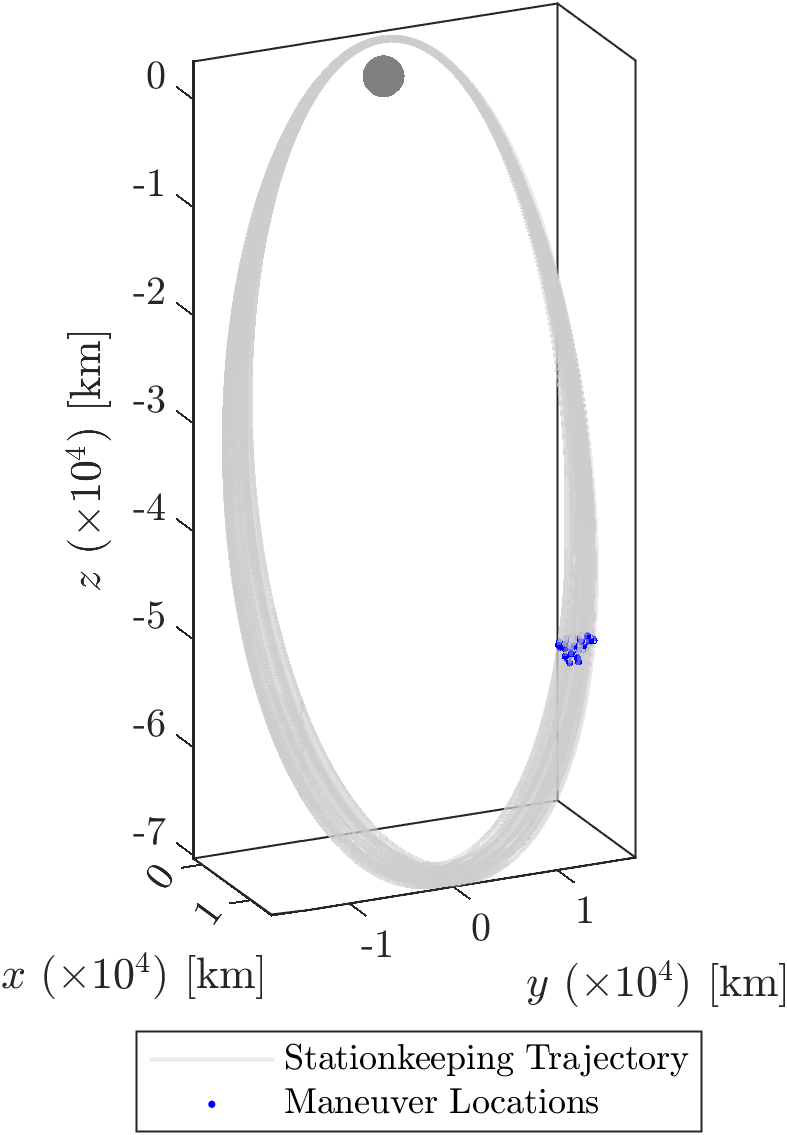}
  \end{minipage}%
  \hfill
  \begin{minipage}{0.5\textwidth}
    \centering
    \begin{subfigure}{\linewidth}
      \includegraphics[width=\linewidth]{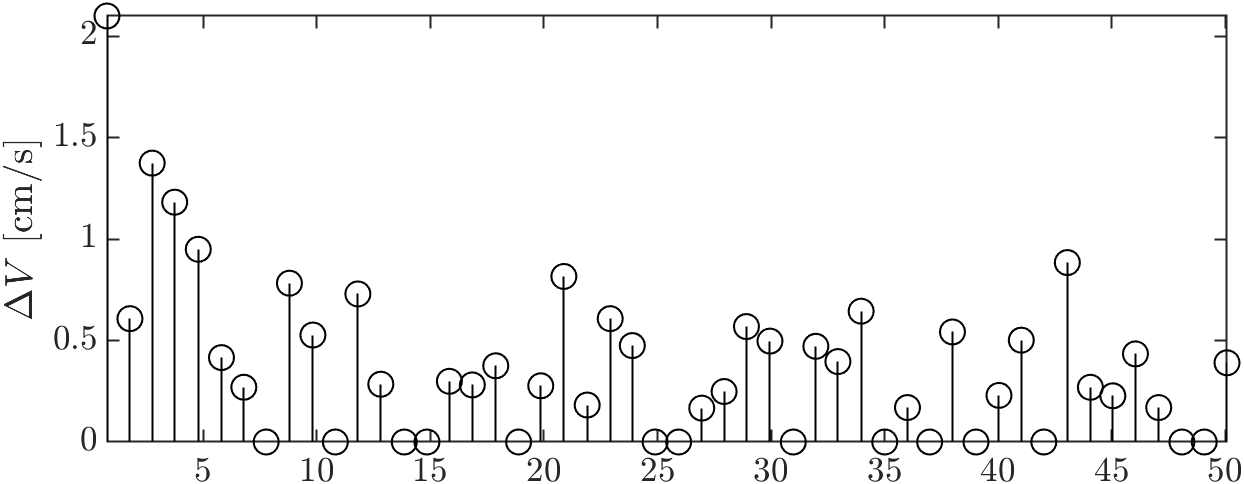}
    \end{subfigure}
    
    \vspace{0cm}
    
    \begin{subfigure}{\linewidth}
      \includegraphics[width=\linewidth]{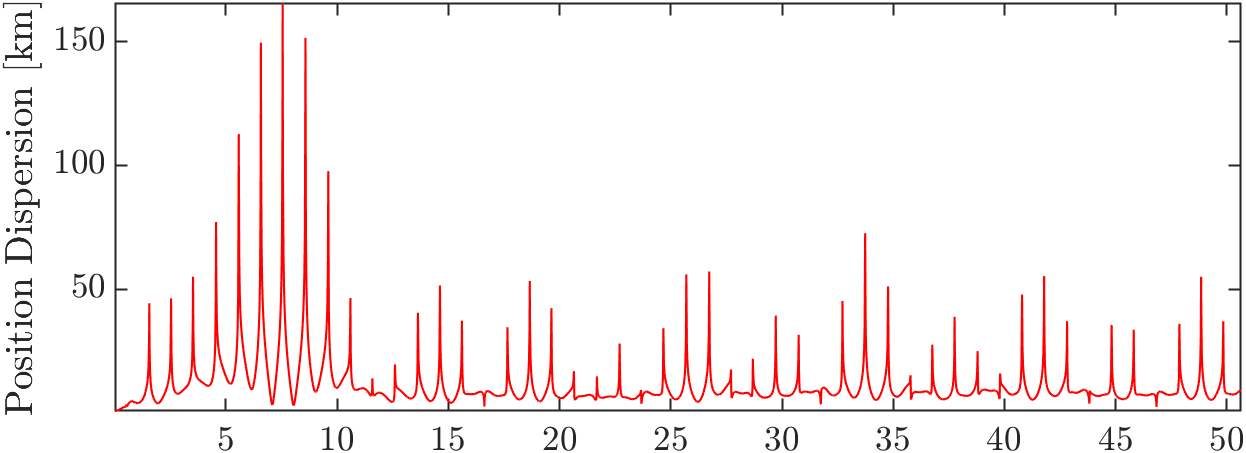}
    \end{subfigure}
    
    \vspace{0cm}
    
    \begin{subfigure}{\linewidth}
      \includegraphics[width=\linewidth]{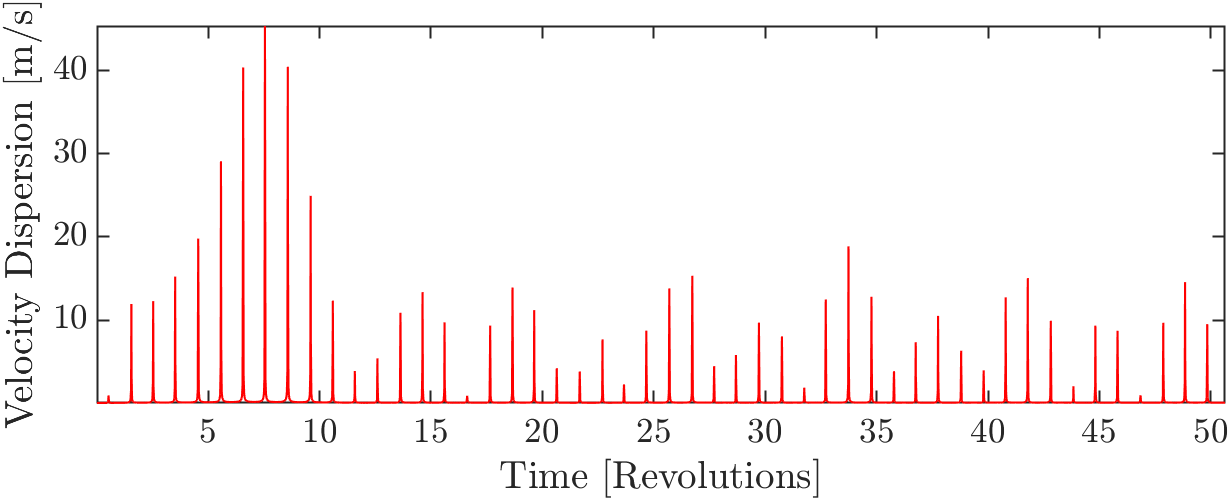}
    \end{subfigure}
  \end{minipage}

  \caption{Absolute control of the Fig.~\ref{fig:Ephem25to6} reference trajectory}
  \label{fig:AbsoluteSK}
\end{figure}

Fig.~\ref{fig:AbsoluteSK} shows the performance of the employed absolute station-keeping scheme about the reference trajectory in Fig.~\ref{fig:Ephem25to6} simulated with a ground truth including solar radiation pressure and with realistic navigation and control errors. Commanded maneuvers with magnitude below 0.15 cm/s are not executed, in line with existing studies of absolute station-keeping \cite{Guzzetti2017,davis2022orbit}.

\subsection{Relative Control}



Given the inadequacies of lower fidelity models in relative control tasks over long durations, relative dynamical motion is generated directly about the reference HFEM trajectory to more accurately represent the relative motion. Since state transition matrices are used to compute the HFEM trajectory and the HFEM QPROs, this relative motion is readily available and could even be stored upon or uplinked to a spacecraft to allow for autonomous relative control. Since dispersions from a reference relative trajectory are typically small in a station-keeping scenario, a simple and effective procedure for computing corrective impulsive maneuvers is to simply target a reference position along the QPRO at some later time \cite{elliott2022phd}. Concretely, to compute a relative maneuver $\Delta V_{rel}$ to deliver the $k$th spacecraft to a desired position $\delta \bm{r}_{f,k}$ (e.g., a reference position along a QPRO) at time $t_f$ from some initial state $\delta \bm{x}_{0,k}$ at $t_0$ simply requires solving
\begin{equation}
    \Delta V_{rel,k} = \Phi^{-1}_{rv}(t_f,t_0) (\delta \bm{r}_{f,k} - \Phi_{rx}(t_f,t_0) \delta \bm{x}_{0,k}).
\end{equation}

\noindent Here, $\Phi_{rv}$ and $\Phi_{rx}$ denote the upper right block and right block of the STM respectively. Under ideal scenarios, the relative state will converge to the reference equilibrium in just two impulses. 
It is efficient to compute, scalable to multiple agents, and independent of the choice of dynamics. Though the control law makes no consideration for fuel expenditure, relative motion is typically so sensitive to changes in velocity in the station-keeping scenario that the commanded $\Delta V$ is usually only on the order of millimeters per second.

\subsection{Safety Concept}

A particular challenge in cislunar formation control is that there exists no known analytical condition which guarantees passive safety of spacecraft in a formation. In Keplerian orbits, passive safety is guaranteed by maintaining perpendicularity between spacecraft relative positions and the direction of the unstable mode. An analytical condition, known as $e-i$ vector separation, has been derived, which when satisfied provides long duration guarantees of safe spacecraft separation under perturbations \cite{damico2006proximity}. In contrast, there does not exist an equally simple and strong notion of passive safety for cislunar orbits. In fact, from a QPRO equilibrium with kilometer level separation, a change in velocity on the order of millimeters per second can deliver a spacecraft to an unsafe configuration in just the span of weeks. In the case in which a spacecraft misses a maneuver, collisions can potentially occur within days or hours. This phenomenon is attributable to the existence of stable and unstable modes of the time-varying eigensystem, whose directions are constantly changing over the span of an orbit.

Collision risk is usually not considered in studies of cislunar formation flight, though solutions that have been proposed include passively-safe sequential convex programming (SCP-DRIFT), whereby nonlinear optimization routines determine a sequence of maneuvers that guarantees passive safety, and maneuver planning leveraging local-toroidal elements (SCP-QPRIT), which allows for elegant suppression of stable and unstable modes \cite{takubo2025safeoptimalnspacecraftswarm}. SCP is computationally expensive and scales poorly with agent count, which makes it impractical for long-duration station-keeping. Additionally, LTCs, which provide weaker guarantees of safety, are not amenable to the inclusion of absolute station-keeping inputs, which are typically spanned by unstable or stable modes.

The safety concept proposed in this work leverages an optimization-based control barrier function (CBF) to provide safety guarantees to the formation. This approach takes advantage of the rarity of collision events, provides recursive safety, and allows for decoupled computation of the absolute and relative control inputs. CBFs typically involve using optimization protocols to ``filter" a nominal control input to be safe \cite{ames2019controlbarrierfunctionstheory}. Though not often applied to spacecraft and typically considered for systems with continuous control, CBFs and similar procedures have been employed for impulsively and continuously controlled spacecraft in Refs.~\citenum{Stephenson2025a} and \citenum{vanwijkrlcbf} for ensuring safety in reinforcement learning routines. The variant of CBF proposed here involves enforcing stochastic keep out zone constraints in a similar manner to passively safe optimal control schemes \cite{guffanti2023passively} and represents one of the first applications of safety filtering to cislunar formation flight.

Let $t_{CO,max}$ and $t_{CO,min}$ be the maximum and minimum cutoff time to a spacecraft's propulsion system following a random missed maneuver event at time $t_0$. It is assumed that following some cutoff time in this interval, the propulsion system will be guaranteed to function at the next available control time. Under these conditions, it is possible to ``filter" the desired control input of all agents at $t_0$ by enforcing keep out zones over the evolution of the agents over the time interval. Particularly, the evolution of the agents subject to commanded control inputs $\{\Delta V_{tot,k} =  \Delta V_{abs} + \Delta V_{rel,k} \}_{k=1}^M$ subject to some initial system state provided by the navigation filter, $\bm{X}_0 \sim \mathcal{N}(\bm{\mu}_0,\Sigma_0 )$, where $\bm{X}_0 = [\bm{x}_1^{\intercal},...,\bm{x}_M^{\intercal} ]^{\intercal}$, $\mathbb{E}(\bm{x}_k) = \bm{\mu_k}$, $\mathbb{E}((\bm{x}_k - \bm{\mu}_k) (\bm{x}_j - \bm{\mu}_j)^{\intercal}) = \Sigma_{kj}$  can be propagated as a distribution with
\begin{equation}
\begin{gathered}
        \bm{\mu}_k(t) = \bm{\mu}_{k,0}  + B \Delta V_k+ \int_{t_0}^{t} \bm{f}(\bm{x}(\tau),\tau) d \tau, \\ 
        \Sigma_{kj}(t) = \Phi_k(t,t_0) (\Sigma_{kj} + \delta_{kj}B G_{exe,k} G^{\intercal}_{exe,j} B^{\intercal}) \Phi^{\intercal}_j(t,t_0) + \delta_{kj} \int_{t_0}^t \Phi_k(t,\tau) G_{rel} G_{rel}^{\intercal} \Phi_k^{\intercal}(t,\tau) d \tau
\end{gathered}
\label{eq:RelProp}
\end{equation}
\noindent Here, $G_{exe,k}$ is an estimate of the maneuver execution error associated with $\Delta V_k$ as determined from Gate's model \cite{Gates1963}, $G_{rel}$ is the covariance associated with the perturbative acceleration to the relative states, $\Phi_k(t,t_0)$ is the STM produced about the mean trajectory of the $k$th agent, $\bm{\mu}_k(t)$, and $\delta_{kj}$ is the Kronecker delta. The above distribution propagation can be simplified by approximating each state mean and covariance as linear with respect to the absolute station-keeping trajectory, whose state and STM is already available. The state can also be propagated in the case of missed maneuvers $\bm{x}_{k,m} \sim \mathcal{N}(\bm{\mu}_{k,m},\Sigma_{k,m})$ simply by setting $\Delta V_k = 0$. To enforce passive safety, it is necessary to ensure that the distributions of each agent are non-overlapping to within some probability for all time until the next available control input. The approach taken here is to enforce that any two spacecraft remain a distance $D_{min}$ apart under $3 \sigma$ certainty \cite{guffanti2023passively}, which can be attained by computing a small modification to the nominal control input through solving a nonconvex optimization problem. Let $t_1$ denote the first control input time after $t_0+t_{CO,max}$ and $t_2$ denote the first control input time after $t_1+t_{CO,max}$. Furthermore, let $\bm{r}_{kj}(t) = R (\bm{x}_k(t) - \bm{x}_j(t))$ be the relative position between spacecraft $k$ and $j$, produced by setting $R= [I_{3\times3}  \ \bm{0}_{3 \times 3} ]$, which is normally distributed according to $\bm{r}_{kj}(t) \sim \mathcal{N}(\bm{p}_{kj}(t),P_{kj}(t))$ where $\bm{p}_{kj}(t) = R ( \bm{\mu}_{k}(t) - \bm{\mu}_j(t) )$ and $P_{kj}(t) = R (\Sigma_k(t) - \Sigma_{kj}(t) - \Sigma_{jk}(t) + \Sigma_j(t) ) R^\intercal$. Note that the covariance of the relative position depends not only the covariances of the states $k$th and $j$th spacecraft, but also the covariance of the states of the states $k$th and $j$th spacecraft with respect to each other. Then, $\{\Delta V_{safe,k}\}_{k=1}^M$ can be solved for from a desired $\{\Delta V_{tot,k}\}_{k=1}^M$ by iteratively solving the problem
\begin{subequations} \label{eq:CBF}
\begin{align}
    \min_{\{\Delta V_{safe,k}\}_{k=1}^M} \quad & \sum_{k=1}^M\| \Delta V_{safe,k} - \Delta V_{tot,k}\|
    \label{eq:CBF_obj} \\
    \text{s.t.} \quad 
    & \bm{p}^{\intercal}_{kj}(t) \bm{p}^{(q-1)}_{kj}(t) - 3 \| P^{\frac{\intercal}{2}}_{kj}(t) \bm{p}^{(q-1)}_{kj}(t) \| \ge D_{min} \| \bm{p}^{(q-1)}_{kj}(t) \|  \ \forall t \in [t_0, t_2], \ \forall k < j \label{eq:CBF_con1} \\
    & \bm{p}^{\intercal}_{kj,m}(t) \bm{p}^{(q-1)}_{kj,m}(t) - 3 \| P^{\frac{\intercal}{2}}_{kj}(t) \bm{p}^{(q-1)}_{kj,m}(t) \| \ge D_{min} \| \bm{p}^{(q-1)}_{kj,m}(t) \|  \ \forall t \in [t_0, t_1], \ \forall k \neq j
    \label{eq:CBF_con2}
\end{align}
\end{subequations}
\noindent Here, $\bm{p}^{(q-1)}_{kj}(t)$ and $\bm{p}^{(q-1)}_{kj,m}(t)$ denote the relative positions at time $t$ as computed from the $(q-1)$th iteration of the scheme. When the scheme is converged, i.e., $\bm{p}^{(q-1)}_{kj}(t) = \bm{p}^{(q)}_{kj}(t)$ and $\bm{p}^{(q-1)}_{kj,m}(t) = \bm{p}^{(q)}_{kj,m}(t)$, the stochastic keep out zones are enforced over the entire duration \cite{guffanti2023passively}. Note that constraint Eq.~\ref{eq:CBF_con1} enforces stochastic keep out zones until the \emph{second} available maneuver time ($t_2$) between all \emph{successfully} maneuvered trajectories whereas constraint Eq.~\ref{eq:CBF_con2} enforces stochastic keep out zones until the \emph{first} available maneuver time ($t_1$) between all possibly \emph{failed} maneuvered trajectories. The former constraint ensures that in the case of two spacecraft missing maneuvers at the next maneuver time ($t_1$), the spacecraft will remain collision free until the following maneuver time ($t_2$), at which time both spacecraft will be capable of maneuvering under the previously made assumption of a maximum propulsion cutoff time. These features of the CBF endow the system with \textit{recursive safety} under the assumption that a spacecraft will be able to execute a maneuver following the cutoff period of an MME.


\begin{algorithm}[t]
\caption{Control Barrier Function Semi-Infinite Program (CBF-SIP)}\label{alg:CBFSIP}
\vspace{-.3cm}
\begin{flushleft}
\hspace*{\algorithmicindent} \textbf{Input} {$\{\bm{\mu}_{0} \}$, $\{ \Sigma_{0}\}$, $t_0,t_1,t_2$, $\{\Delta V_{tot,k}\}$} \\
\hspace*{\algorithmicindent} \textbf{Output} $\{\Delta V_{safe,k}\}$
\end{flushleft}
\vspace{-.4cm}
\begin{algorithmic}[1]
\Function{\texttt{CBF-SIP}}{}
    \State $\{\Delta V_{safe,k}\} \gets \{\Delta V_{tot,k}\}$, $\{ t_c \} \gets \emptyset$
    \While{Maximum Iterations not Exceeded}
        \State $\{ \{\bm{\mu}_{k}(t) \}, \{\bm{\mu}_{k,m}(t) \}, \{\Sigma_{kj}(t) \} \} \gets $ Propagate \ref{eq:RelProp} with $\{\Delta V_{safe,k}\}$.
        \State \textbf{if} \ref{eq:CBF_con1} and \ref{eq:CBF_con2} are satisfied abouth current reference \textbf{then return} $\{\Delta V_{safe,k}\}$.
        \State $\{ t_c\} \gets \{ t_c\} \cup $ Time points where \ref{eq:CBF_con1} or \ref{eq:CBF_con2} are violated and their LHS are locally minimal.
        \State $\{\Delta V_{safe,k}\}  \gets $ Solve \ref{eq:CBF} with linearized version of \ref{eq:CBF_con1} and \ref{eq:CBF_con2} $\forall t \in \{ t_c\}$.
    \EndWhile
\EndFunction
\end{algorithmic}
\end{algorithm}

The optimization problem in Eq.~\ref{eq:CBF} is a non-convex semi-infinite problem (SIP). To solve it efficiently, an iterative approach that is commonly used to solve SIPs is employed whereby affine inequality constraints are only applied at discrete time points at which the semi-infinite constraint is violated \cite{Shapiro01022009}. These time points are selected by finding the times at which the worst violation of the inequality constraints occurs. This recursive procedure is outlined in Alg.~\ref{alg:CBFSIP}. In practice, the authors do not find it necessary to augment the non-convex constraints with slack or introduce a trust region for robust convergence. The CBF-SIP should be \textit{recursively feasible} so long as the covariance at the next maneuver time does not become larger than the predicted covariance from the previous maneuver time and the estimated relative state error does not change drastically from one iteration to the next.

Typically, a maneuver profile is already free of collision risk, in which case the CBF returns the same input profile without running any optimization routine. This can be contrasted to the SCP-DRIFT and SCP-QPRIT techniques, which do not exploit the rarity of collision events for efficient solution procedures. Moreover, leveraging the SIP reduction technique means that each optimization problem involves very few constraints, thereby improving scalability with respect to the number of agents. The CBF-SIP, however, does not necessarily provide guarantees of stability when combined with absolute and relative maneuver algorithms. In the case of very tight formations, where collision risk is higher, the CBF-SIP may result in greater dispersions and higher $\Delta V$ expenditure. For this reason, the CBF-SIP can be thought of as an additional form of (safe) disturbance to the formation that has magnitude proportional to the tightness of the formation.

\subsection{Maneuver Time Placement}

The performance of the absolute control, relative control, and safety filter are heavily dependent on the placement of maneuver times. The nature of NRHO trajectories dictates that maneuvers be performed near apoapsis crossing, where the dynamics are far less sensitive so as to not amplify navigation and control execution errors \cite{davis2022orbit,shimane2025output}. Applying maneuvers within the time intervals corresponding to true anomalies between $160\deg$ and $200\deg$, a duration of $\sim3.6$ days, provides ample maneuver freedom without incurring large errors due to dynamical sensitivity. This duration is also sufficiently long to allow for the resolution of and recovery from MMEs by ground operators for absolute control maneuvers, which have the greatest magnitude and therefore the most severe effect when missed. It is also best for absolute maneuvers, which are applied once per orbit, to be performed as early as possible in this interval, so that the ensuing relative maneuvers can be used to clean up formation dispersions introduced by the absolute maneuvers. This ensures that the magnitude of the relative dispersions at perilune crossing, the most sensitive and chaotic region, is low. Thus, absolute maneuvers are performed at times corresponding to a true anomaly of $160\deg$ and several relative maneuvers are performed up to and including a true anomaly of $200\deg$. The number of relative maneuvers is a design parameter dependent on desired $\Delta V$ expenditure, which increases with maneuver frequency, and tracking performance, which improves with maneuver frequency.

\subsection{Virtual Chief Placement}

The choice of absolute state to represent the formation can have a dramatic effect on the $\Delta V$ consumption of the formation and the extent to which this consumption is spread among the agents. For example, a virtual chief can be selected, for which the deputy spacecraft maintain a relative orbit with respect to, such that the overall fuel expenditure or distribution of fuel expenditures is minimized. This has been demonstrated in circular orbits \cite{LIPPE2021162} and eccentric orbits \cite{lippevcplacement} for Keplerian motion. Especially for large swarms, it may be detrimental to have a specific spacecraft operate as the chief, since this relocates the disturbances on the chief to the deputies, removing relative station-keeping burden from the chief and placing it on the deputies. In general, the $\Delta V$ associated with the choice of virtual chief is difficult to determine since both absolute and relative $\Delta V$ expenditure is affected. However, a simple and effective strategy that distributes relative station-keeping burden among the formation is to place the virtual chief at the mean of the state estimates. If the collection of absolute state estimates from the navigation system are given by $\{ \hat{\bm{x}}_k \}_{k=1}^M$, then the virtual chief state and deputy states are determined by
\begin{equation}
    \bm{x}_c(t) = \frac{1}{|\mathcal{I}(t)|} \sum_{k \in \mathcal{I}(t)} \hat{\bm{x}}_k(t) - \delta \bm{x}_{ref,k}(t)
\end{equation}
\noindent where $\delta \bm{x}_{ref,k}(t)$ is the relative reference state of deputy $k$ at time $t$ along the QPRO and $\mathcal{I}$ is an index set over which the averaging is applied. The index set $\mathcal{I}$ at time $t$ includes the indices for the agents for which a maneuver was available at the last control input (i.e., no propulsion cutoff). This is done so that the occurrence of a missed maneuver from any agent does not drastically alter the location of the virtual chief, thereby allocating relative station-keeping burden to the agents that missed their last maneuver. Relative state estimates are computed with respect to the virtual chief with
\begin{equation}
    \delta \hat{\bm{x}}_k = \hat{\bm{x}}_k(t) - \bm{x}_c(t).
\end{equation}

This averaging procedure effectively distributes the disturbance and uncertainty among all agents. If the absolute state estimate covariances are approximately equal, then the relative state uncertainty should be approximately equivalent for each agent. Moreover, the effect of disturbance and control execution error before the current time should result in approximately the same relative dispersion with this choice of virtual chief. 

\section{Navigation Simulation}

\begin{wraptable}{r}{8cm}
{\small
\centering
\begin{tabular}{cc}
\Xhline{0.4pt}
\Xhline{0.4pt}
Parameter & Value ($3\sigma$)   \\ 
\hline
Initial Absolute Position Navigation Error & 1 [km]   \\ 
Initial Absolute Velocity Navigation Error & 1 [cm/s]  \\ 
Initial Relative Position Navigation Error & 1 [m]  \\ 
Initial Relative Velocity Navigation Error & 0.1 [mm/s]  \\ 
Absolute Pseudo-range (DSN) Error  & 1 [m]  \\ 
Absolute Pseudo-range Rate (DSN) Error & 0.1 [mm/s]   \\ 
Relative Pseudo-range (IS) Error  & 1 [m]  \\ 
Relative Pseudo-range Rate (IS) Error & 0.1 [mm/s]   \\
\hline
\hline
\end{tabular}
\caption{Navigation System Parameters}
\label{tab:navigation}
}
\end{wraptable}

To accurately represent and simulate the performance of the station-keeping algorithms, realistic navigation algorithms are implemented. Navigation updates are given in the form of absolute and relative pseudo-range and pseudo-range rates, representing Deep Space Network (DSN) and intersatellite (IS) measurements and fed to an Extended Kalman Filter (EKF). The measurement architectures explored include a cluster architecture, whereby measurements include the absolute range of a hub spacecraft and the relative ranges between all spacecraft, a fully connected architecture, whereby measurements include the absolute ranges of all spacecraft and relative ranges between all spacecraft, and a hub-spoke architecture, whereby measurements include absolute range of a hub spacecraft and relative range between the hub and spoke spacecraft. Absolute ranging is scheduled to occur at -72 hours, -48 hours, -7 hours and +12 hours with respect to an absolute maneuver time (true anomaly of $160\deg$), with 10 measurements taken over the span of one hour, in line with the measurement frequency proposed in Ref.~\citenum{shimane2025output}. In addition, relative ranging is scheduled for the same time as absolute ranging and 1 hour before any relative maneuver time for the fully connected architecture, every 24 hours for the cluster architecture, and every 6 hours for the hub and spoke architecture. Tab.~\ref{tab:navigation} contains the relevant navigation parameters applied to all ensuing test cases.




\section{Test Scenarios and Results}

To demonstrate the flexibility and scalability of the proposed station-keeping architecture to relevant cislunar formation flight scenarios, the algorithms are deployed on and simulated in the following test scenarios:

\begin{enumerate}[itemsep=0pt, topsep=0pt]
    \item[1.] A two spacecraft fully connected formation with $\sim10$ km separation under a noisy configuration.
    \item[2.] A five spacecraft cluster formation with $\sim1$ km relative separation under a quiet configuration.
    \item[3.] A nine spacecraft hub-spoke formation with $\sim100$ km relative separation under a noisy configuration.
\end{enumerate}

\noindent The parameters defining a quiet and noisy configuration are given in Tab.~\ref{tab:SimParams}. 

\begin{table}[t]
{\small
\centering
\begin{tabular}{ccc}
\Xhline{0.4pt}
\Xhline{0.4pt}
Parameter & Quiet Configuration & Noisy Configuration  \\ 
\hline
Fixed Magnitude Execution Error & 0.01 [mm/s] & 0.3 [mm/s]  \\ 
Proportional Magnitude Execution Error & 0.1 [\%] & 1 [\%]  \\ 
Fixed Pointing Execution Error & .01 [mm/s] & 0.3 [mm/s]  \\ 
Proportional Pointing Execution Error & 0.001 [rad] & 0.01 [rad]  \\ 
Missed Maneuver Event Rate  & 1 [\%] & 5 [\%]  \\ 
Maximum Propulsion Cutoff Time & 3 days & 3 days \\ 
Minimum Propulsion Cutoff Time & 1 days & 1 days \\ 
Minimum Impulse Bit & 0.01 [mm/s] & 0.3 [mm/s] \\
Solar Radiation Acceleration & 1e-5  [mm/$\text{s}^2$]  $\pm$ 5 \% & 5e-5  [mm/$\text{s}^2$]  $\pm$ 10 \% \\
Minimum Safe Separation & 100 [m] & 1000 [m] \\
\hline
\hline
\end{tabular}
\caption{Control and simulation parameters}
\label{tab:SimParams}
}
\end{table}
\begin{table}[t]
{\footnotesize
\centering
\setlength{\tabcolsep}{1pt}
\begin{tabular}{cccc}
\Xhline{0.4pt}
\Xhline{0.4pt}
Metric & Test Case 1, Fig.~\ref{fig:TwoAgents} & Test Case 2, Fig.~\ref{fig:FiveAgents} & Test Case 3, Fig.~\ref{fig:NineAgents}    \\ 
\hline
Monte Carlo Trials & 10 & 10 & 10 \\
Computed $\Delta V_{abs}$, [cm/s/year] & $53.99 \pm 7.77$ & $9.08\pm2.61$ & $52.63\pm 5.44$ \\
\makecell{Total Executed $\Delta V$ \\ per agent, [cm/s/year]} &
{\footnotesize $
\left[
\begin{array}{@{}l@{}}
81.39 \pm 9.42 \\
75.34\pm 9.49
\end{array}
\right]
$
}& 
{\footnotesize $
\left[
\begin{array}{@{}l@{}}
29.65 \pm 7.47, 23.63 \pm 5.62 \\
24.76 \pm 4.33, 23.85 \pm 5.91 \\
24.46 \pm 5.64
\end{array}
\right]
$
}
 & {\footnotesize $
\left[
\begin{array}{@{}l@{}}
75.50 \pm 5.30, 82.76 \pm 7.15 \\
88.03 \pm 10.17, 89.99 \pm 22.93 \\
96.85 \pm 17.91, 89.43 \pm 8.76 \\
93.08 \pm 11.58, 99.58\pm10.78 \\
100.29 \pm 7.74
\end{array}
\right]
$
}
\\
Relative Separation Ranges [km] & [1.04,57.90] & [0.099,7.34] & [2.38,523.98] \\
CBF-SIP Activation Rate [\%] & 2.8 & 29.13 & 0.01 \\
\makecell{Absolute Maneuver \\ Computation Time [sec]} & $63.40 \pm 29.65$ & $45.34\pm18.89$ & $84.64\pm43.03$ \\
\makecell{Relative Maneuver \\ Computation Time [millisec]} & $0.49\pm0.30$ & $0.60\pm0.77$ & $1.44\pm4.92$ \\
CBF-SIP Computation Time [sec] & $0.07\pm0.32$ & $0.89\pm1.36$ & $0.89\pm0.58$ \\
\hline
\hline
\end{tabular}
\setlength{\tabcolsep}{6pt}
\caption{Performance metrics of each test case with standard deviations}
\label{tab:metrics}
}
\end{table}

Simulations are performed using the NRHO reference from Fig.~\ref{fig:Ephem25to6} for 50 revolutions starting at an epoch of 2024-10-29 12:00:00 TDB. The simulation process is outlined as follows:

\begin{enumerate}[itemsep=0pt, topsep=0pt]
    \item Uniformly sample a reflective coefficient associated with each agent in accordance with the SRP errors in Tab.~\ref{tab:SimParams}. Propagate the agents until the next maneuver time. Propagate the navigation filter with measurements applied at relevant times.
    \item Compute the estimated state of the virtual chief as the mean state estimate of all agents. Compute estimated relative states with respect to virtual chief.
    \item If the maneuver time is an absolute maneuver time, compute an absolute maneuver for the virtual chief. 
    \item Compute the relative maneuver for each agent as described.
    \item Sum the absolute and relative maneuvers and apply the CBF-SIP safety filter in Alg.~\ref{alg:CBFSIP}.
    \item Remove maneuvers that violate the minimum impulse bit defined in Tab.~\ref{tab:SimParams}. Apply maneuver execution errors from Gate's model \cite{Gates1963} with errors from Tab.~\ref{tab:SimParams}. If a uniformly sampled random number is below the MME rate, then uniformly sample a propulsion cutoff time between the minimum and maximum bounds. All maneuvers for the given agent are set to zero for the remainder of this cutoff time. Add the resulting $\Delta V$ to the state of each agent.
\end{enumerate}

\begin{figure}[t]
  \centering
  
  \begin{subfigure}{0.59\textwidth}
    \centering
    \begin{subfigure}{\linewidth}
      \centering
      \includegraphics[width=\linewidth]{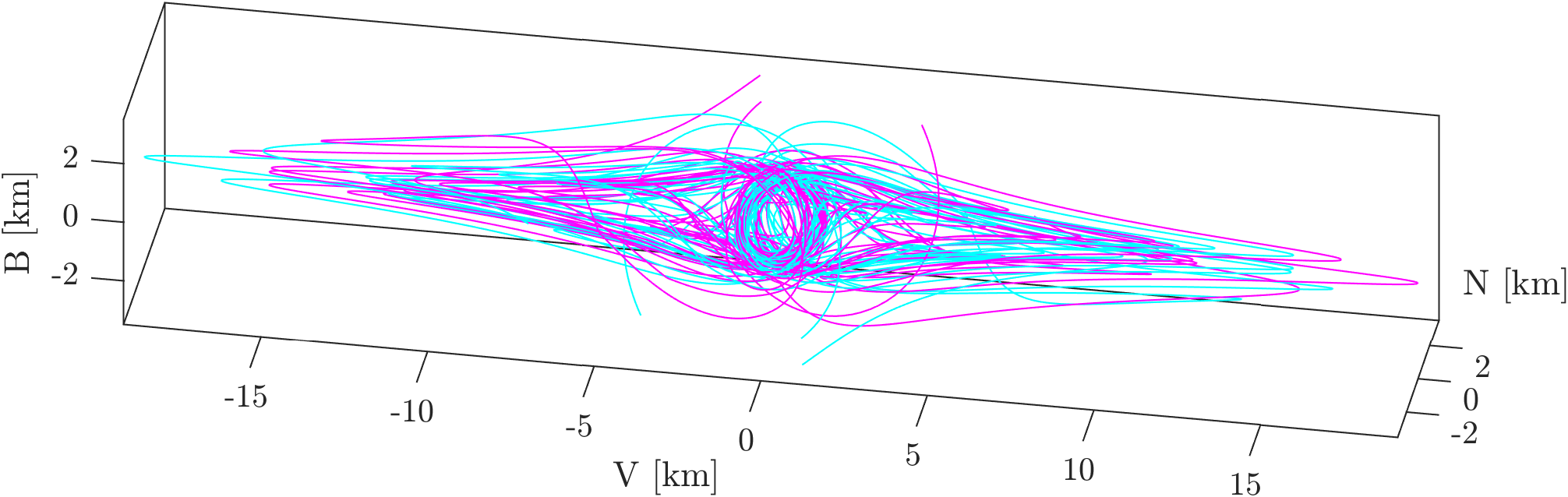}
      \vspace{-.4cm}
    \end{subfigure}
    \begin{subfigure}{.58\linewidth}
      \centering
      \includegraphics[width=\linewidth]{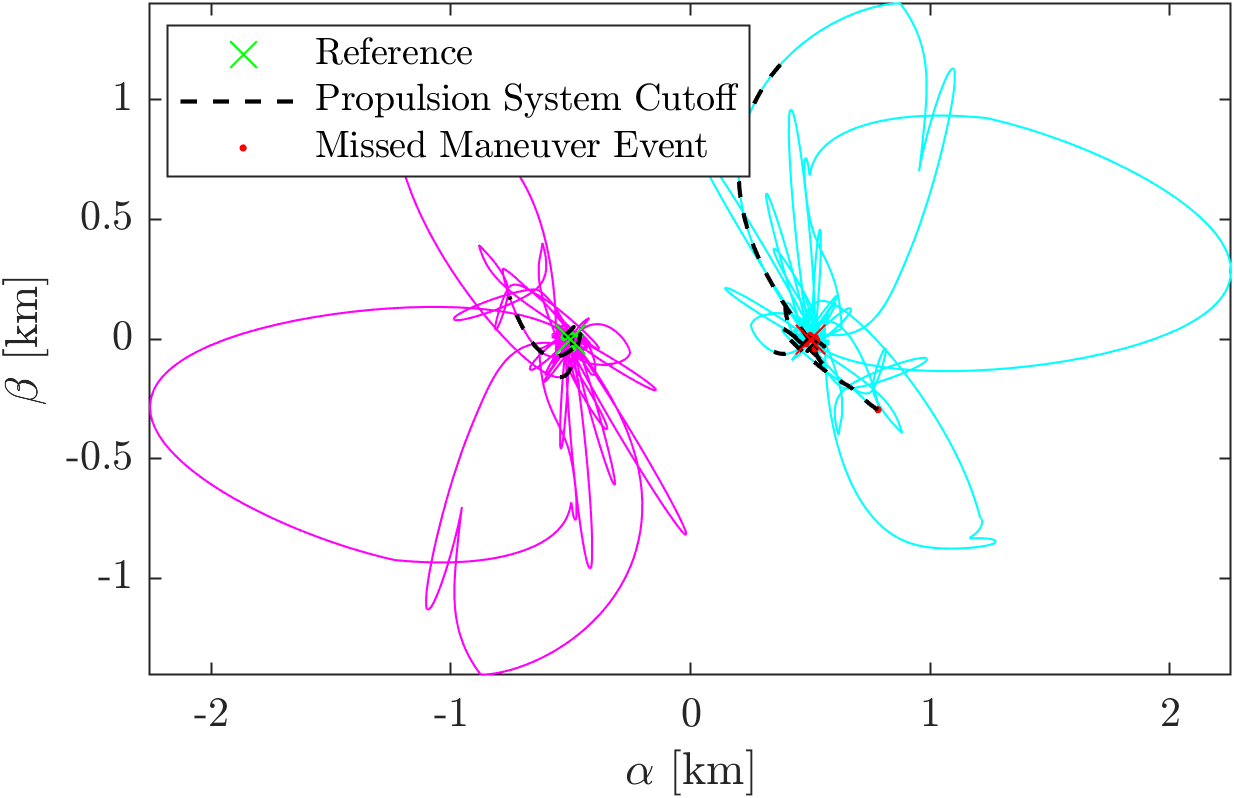}
    \end{subfigure}
    \hfill
    \begin{subfigure}{.39\linewidth}
      \centering
      \includegraphics[width=\linewidth]{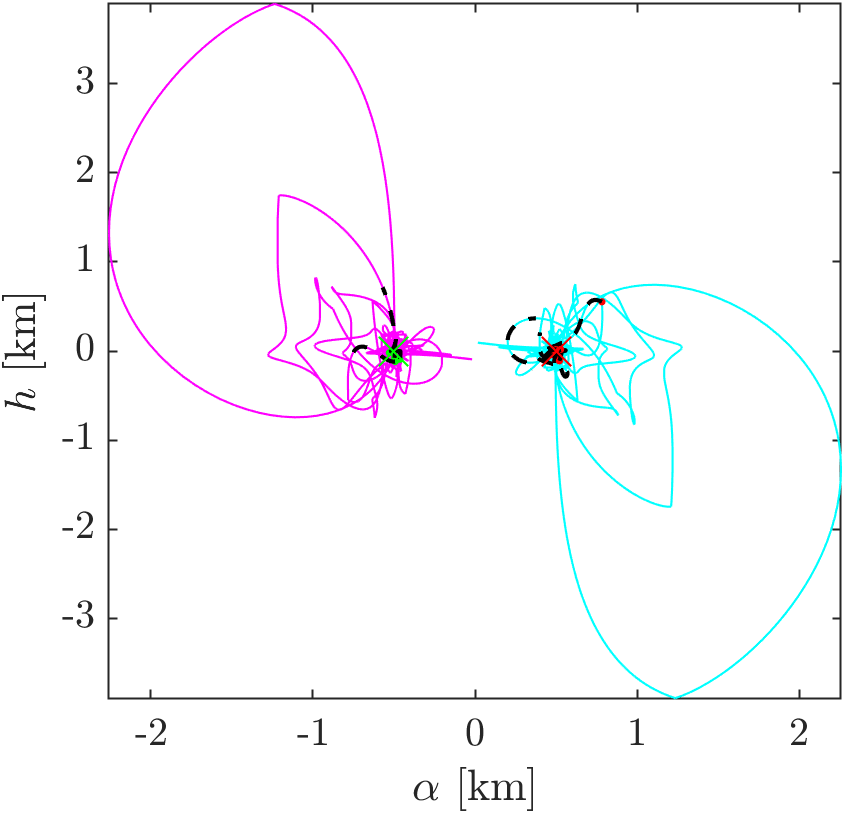}
    \end{subfigure}
    \caption{Relative Positions and Local Toroidal Coordinates}
  \end{subfigure}
  \hfill
  \begin{subfigure}{0.4\textwidth}
    \centering
    \includegraphics[width=\linewidth]{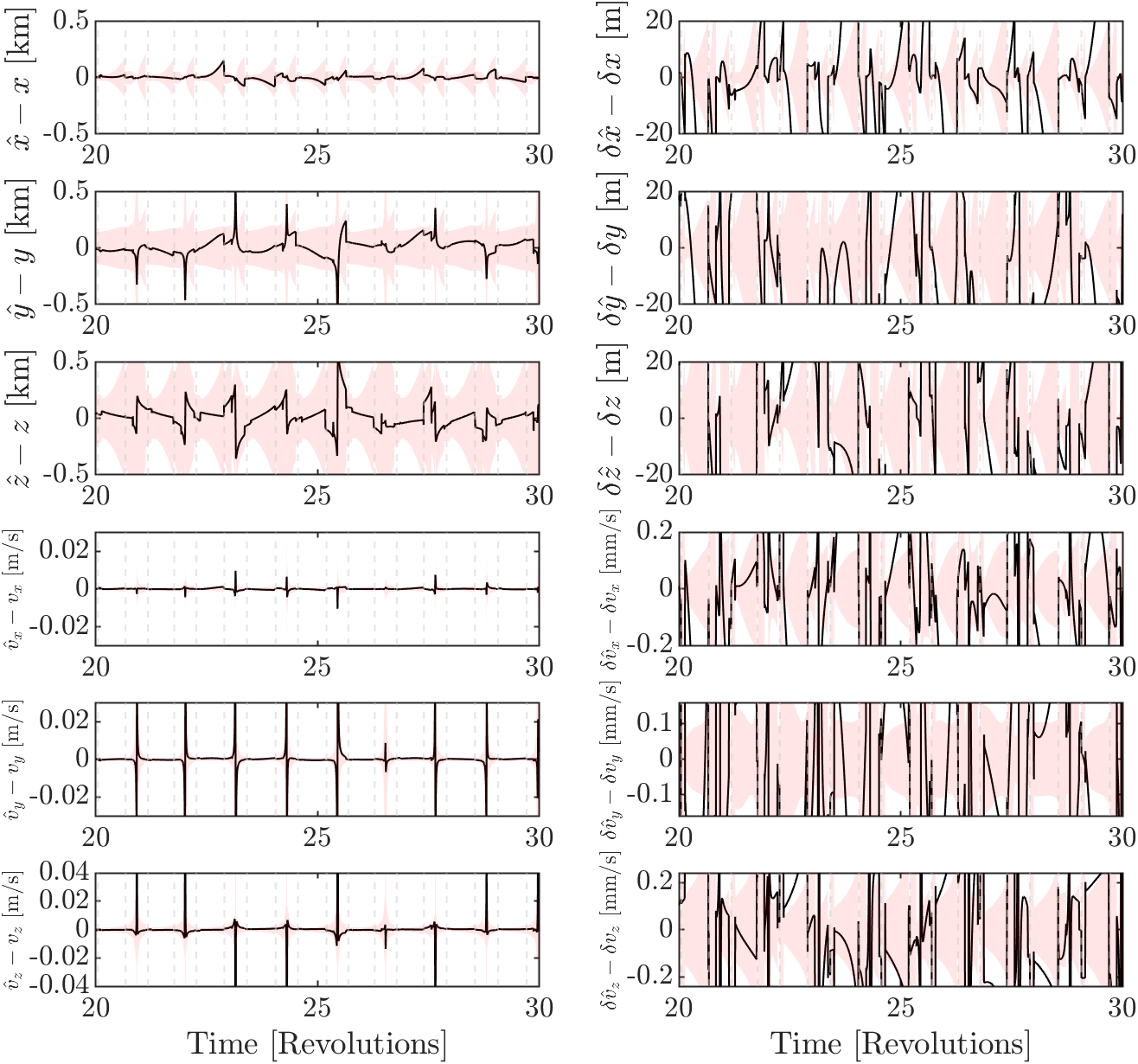}
    \caption{Navigation Subset}
  \end{subfigure}
  
  \vspace{0em}

  \begin{subfigure}{0.3\textwidth}
    \centering
    \includegraphics[width=\linewidth]{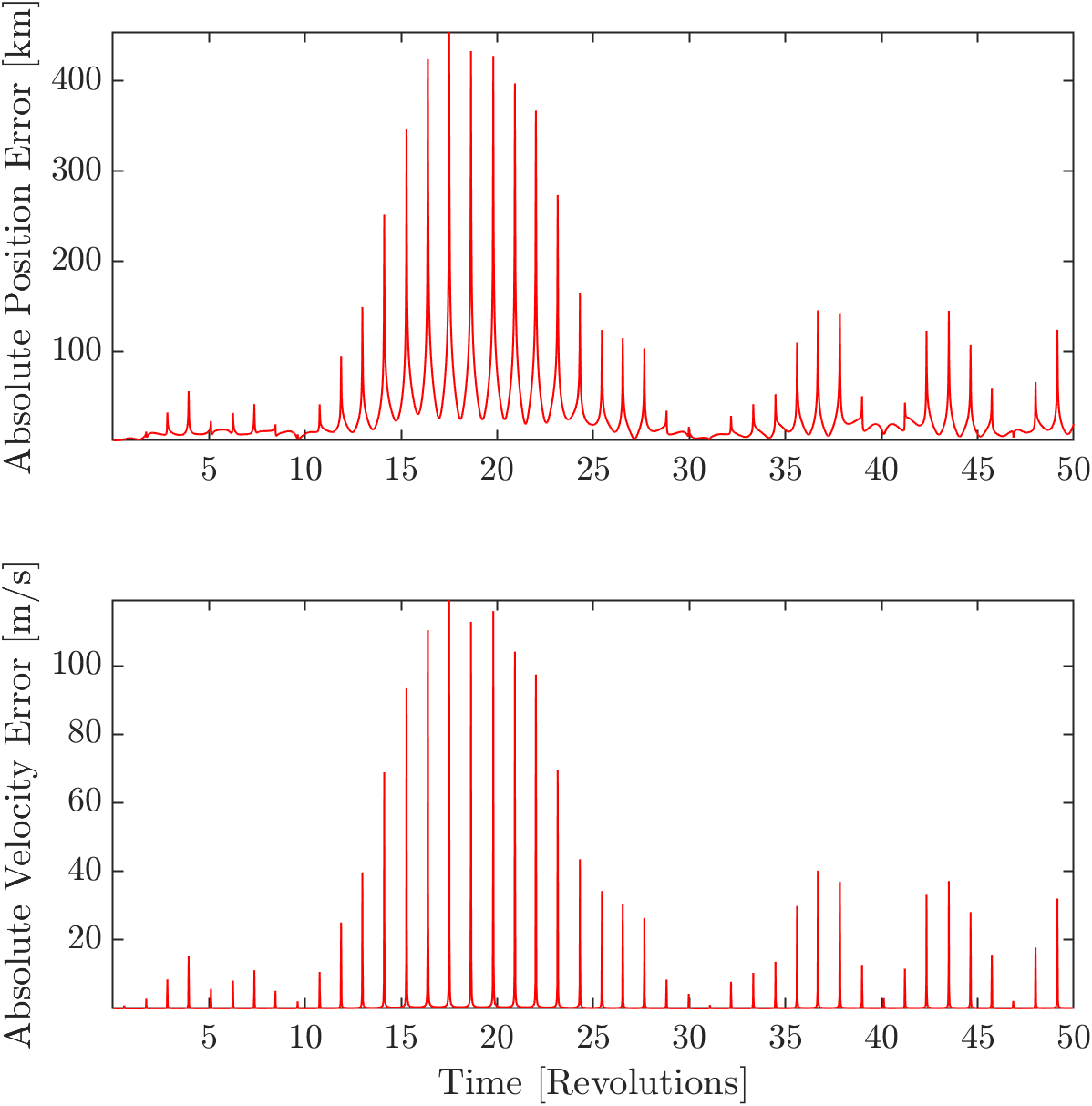}
    \caption{Absolute Reference Tracking Performance}
  \end{subfigure}
  \hfill
  \begin{subfigure}{0.3\textwidth}
    \centering
    \includegraphics[width=\linewidth]{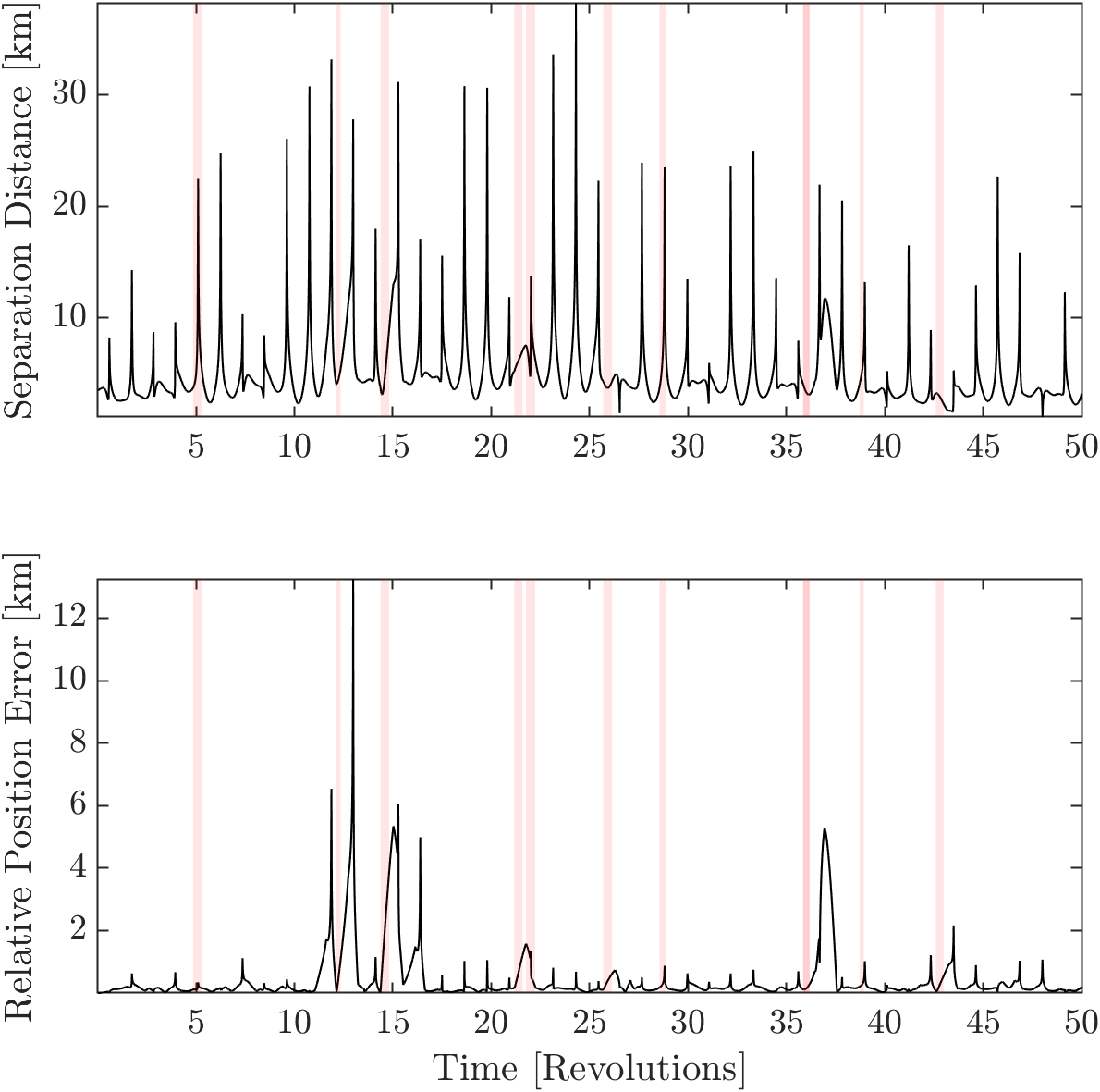}
    \caption{Relative Reference Tracking and Separation}
  \end{subfigure}
  \hfill
  \begin{subfigure}{0.32\textwidth}
    \centering
    \begin{subfigure}{.98\linewidth}
      \centering
      \includegraphics[width=\linewidth]{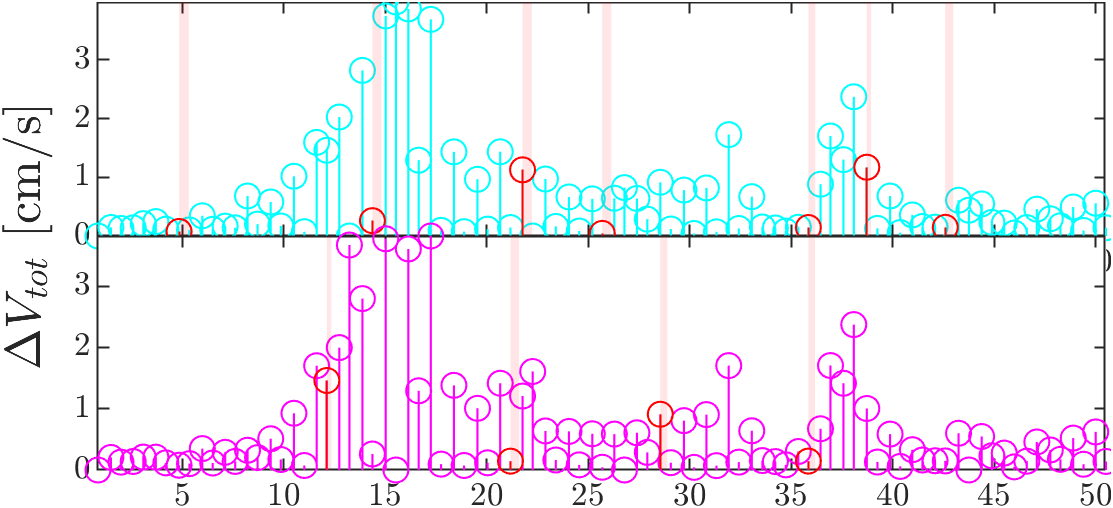}
    \end{subfigure}
    \begin{subfigure}{.98\linewidth}
      \centering
      \includegraphics[width=\linewidth]{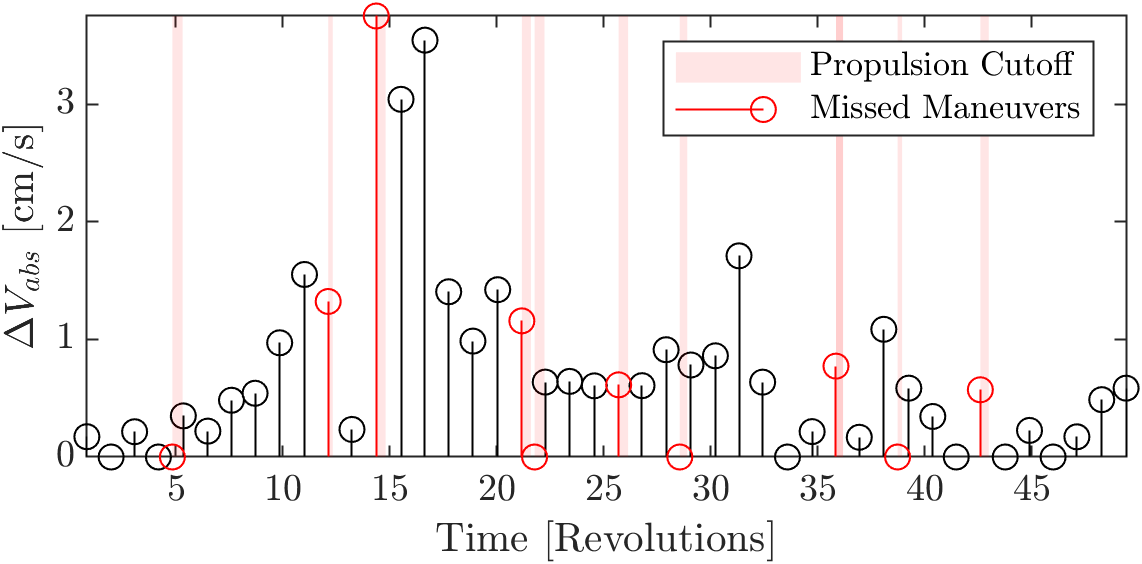}
    \end{subfigure}
    \caption{Absolute and Total $\Delta V$ Expenditure}
  \end{subfigure}
  
  \caption{Test Case 1 Trial: Two agent station-keeping for 50 revolutions}
  \label{fig:TwoAgents}
  \vspace{-.7cm}
\end{figure}

\begin{figure}[t]
  \centering
  
  \begin{subfigure}{0.55\textwidth}
    \centering
    \begin{subfigure}{1\linewidth}
      \centering
      \includegraphics[width=\linewidth]{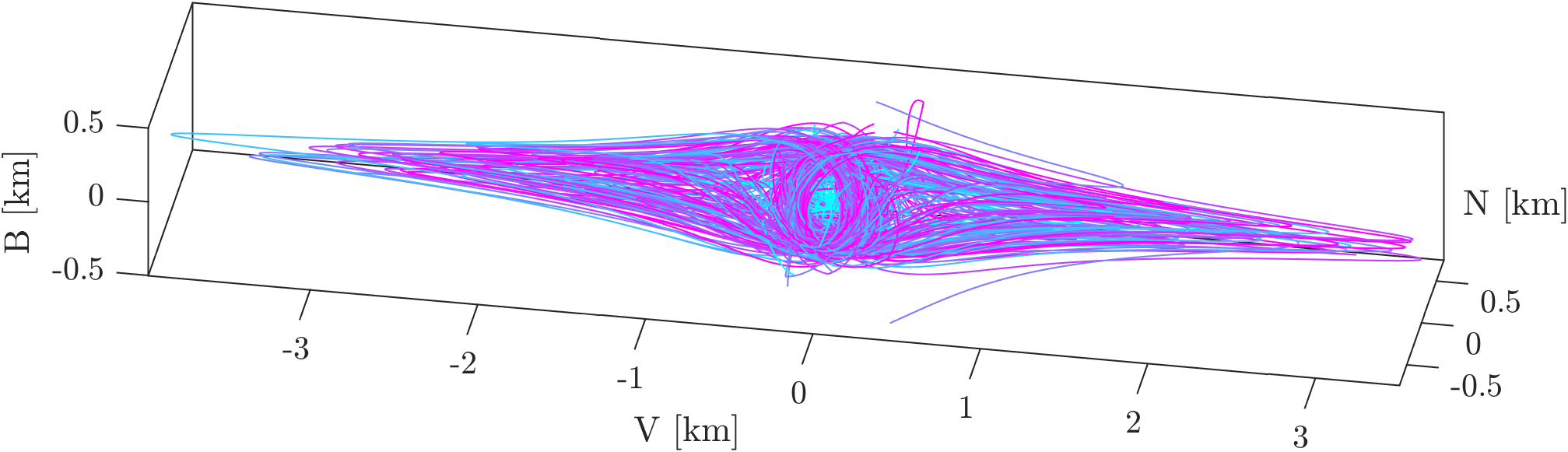}
    \end{subfigure}
    \vspace{0.2em}
    \begin{subfigure}{.61\linewidth}
      \centering
      \includegraphics[width=\linewidth]{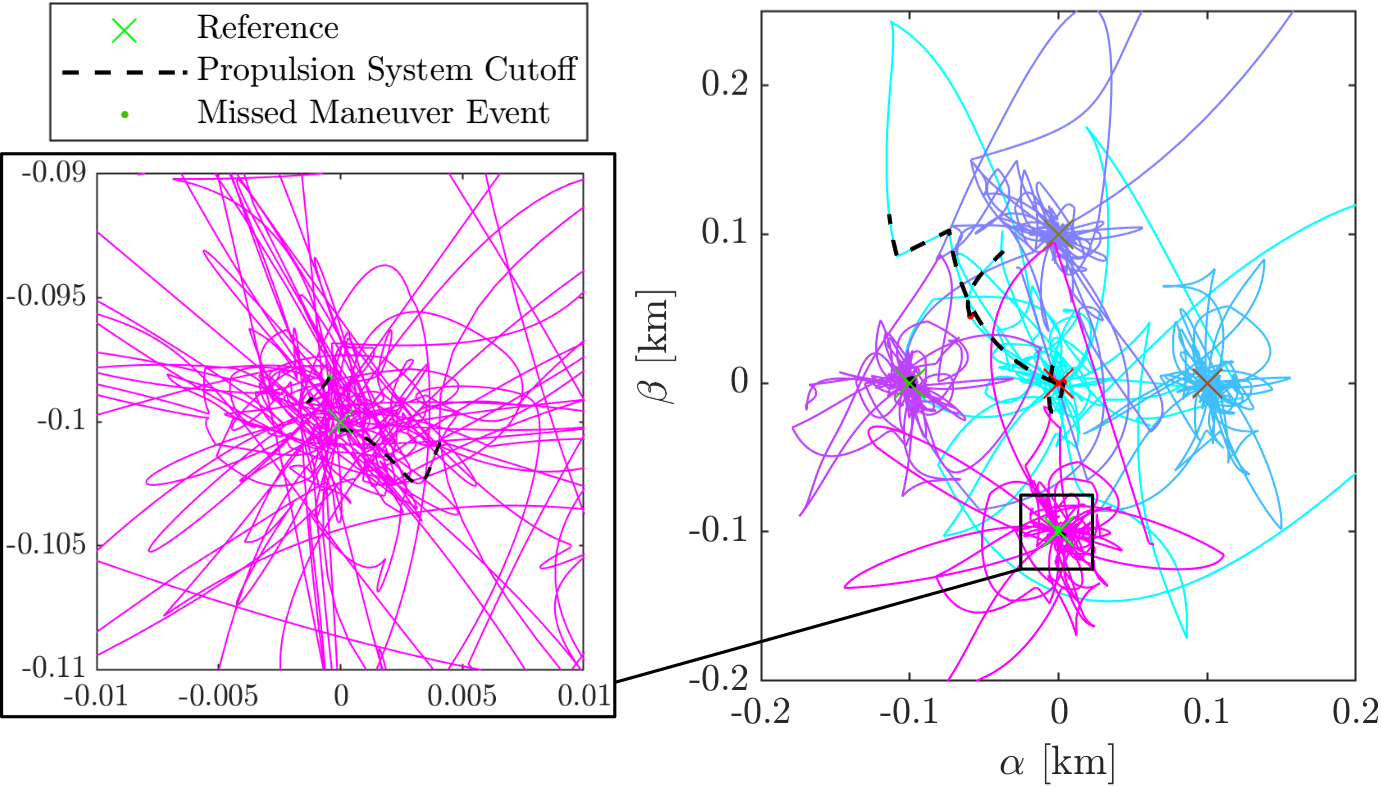}
    \end{subfigure}
    \begin{subfigure}{.3\linewidth}
      \centering
      \includegraphics[width=\linewidth]{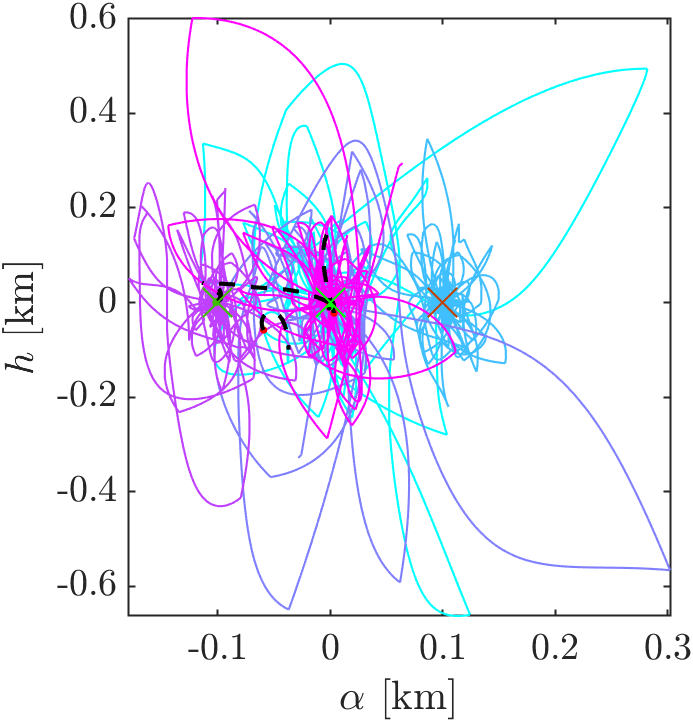}
    \end{subfigure}
    \caption{Relative Positions and Local Toroidal Coordinates}
  \end{subfigure}
  \hfill
  \begin{subfigure}{0.4\textwidth}
    \centering
    \includegraphics[width=\linewidth]{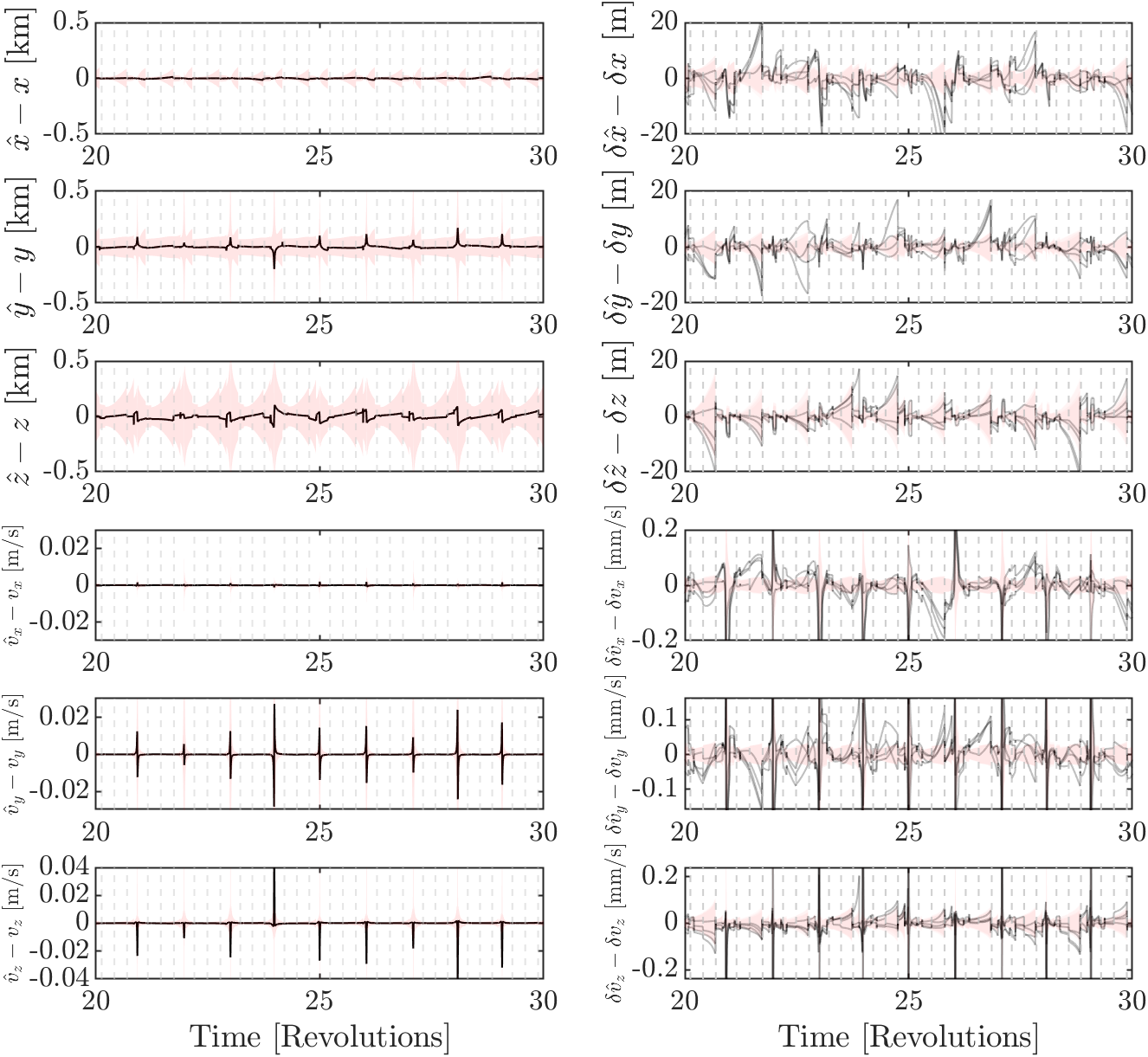}
    \caption{Navigation Subset}
  \end{subfigure}
  
  \vspace{0em}

  \begin{subfigure}{0.3\textwidth}
    \centering
    \includegraphics[width=\linewidth]{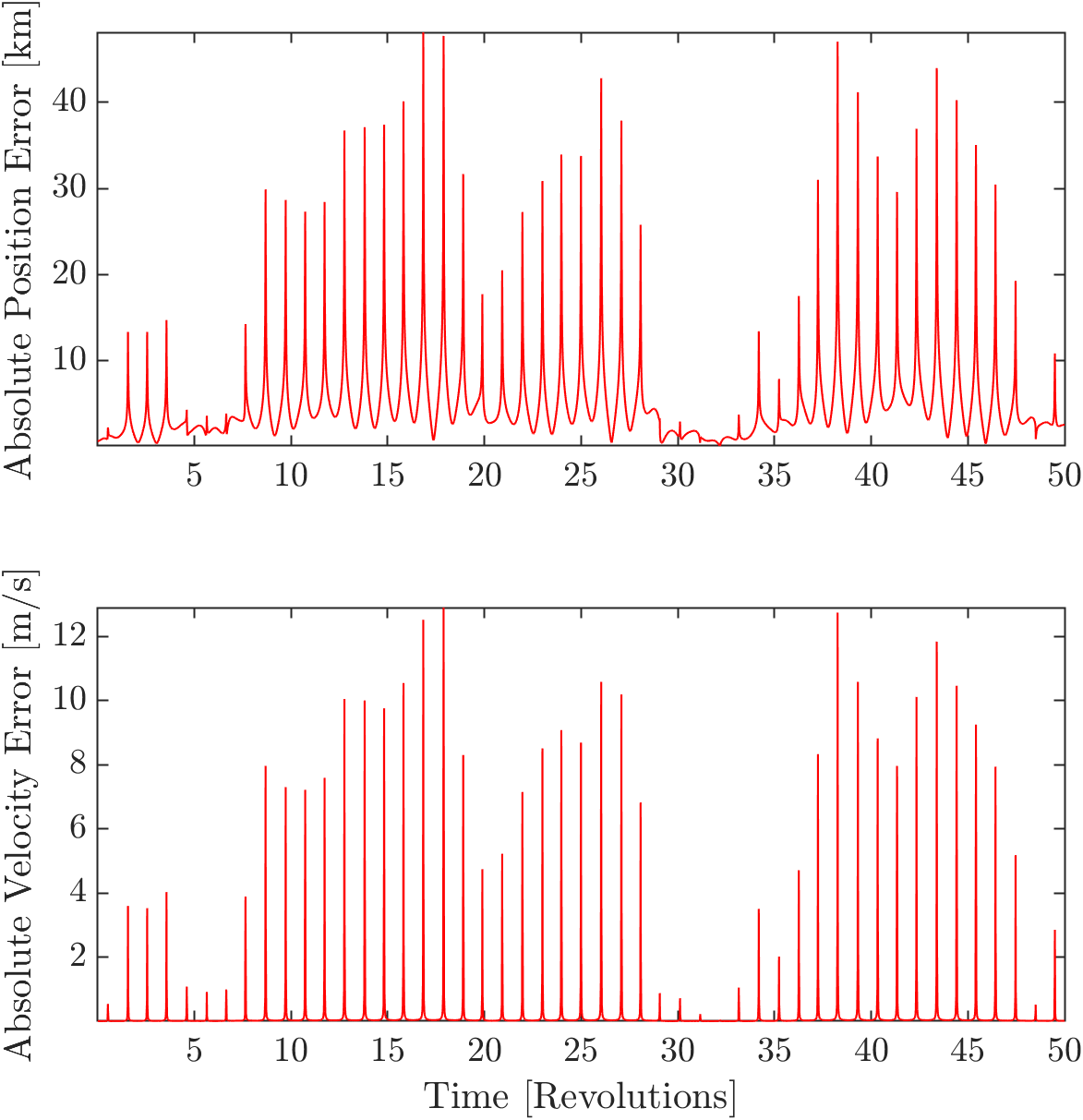}
    \caption{Absolute Reference Tracking Performance}
  \end{subfigure}
  \hfill
  \begin{subfigure}{0.3\textwidth}
    \centering
    \includegraphics[width=\linewidth]{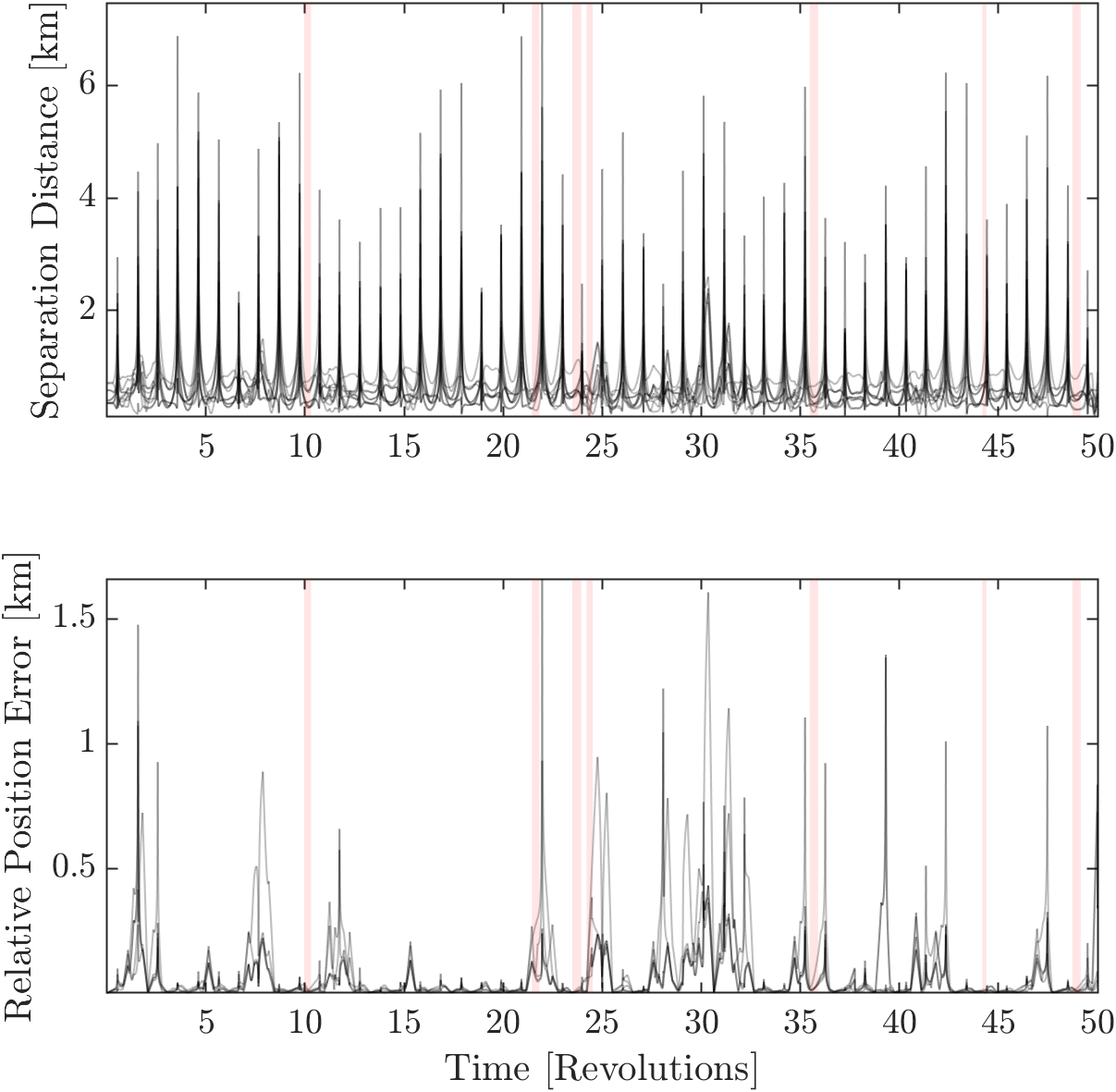}
    \caption{Relative Reference Tracking and Separation}
  \end{subfigure}
  \hfill
  \begin{subfigure}{0.32\textwidth}
    \centering
    \begin{subfigure}{1\linewidth}
      \centering
      \includegraphics[width=\linewidth]{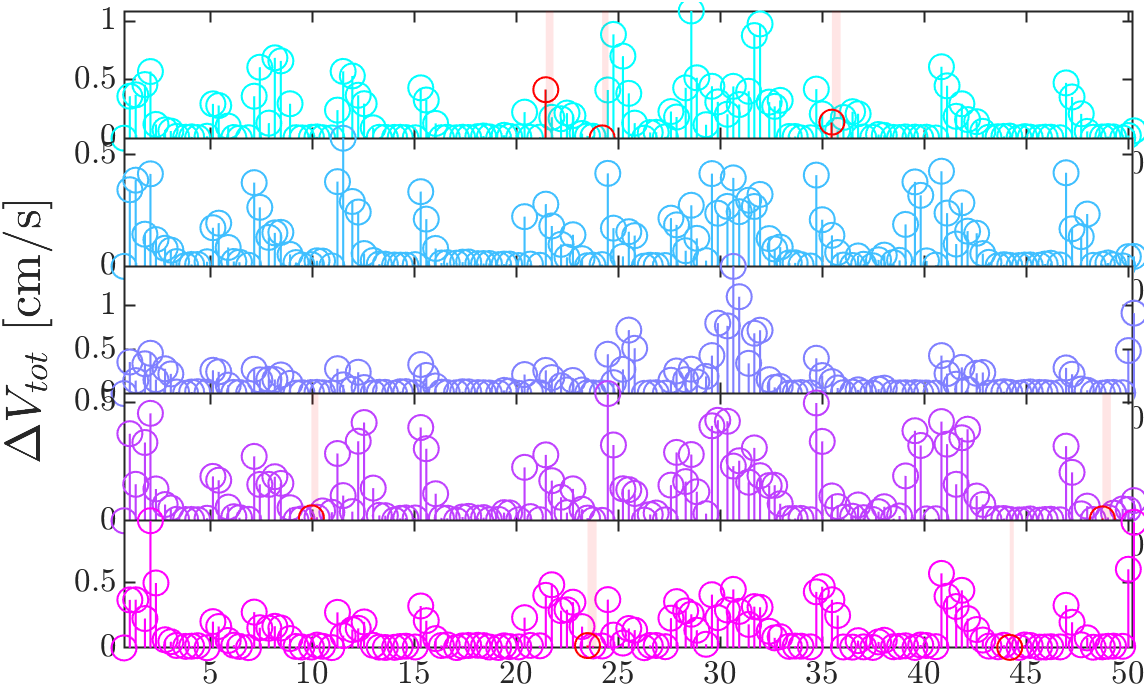}
    \end{subfigure}
    \begin{subfigure}{\linewidth}
      \centering
      \includegraphics[width=\linewidth]{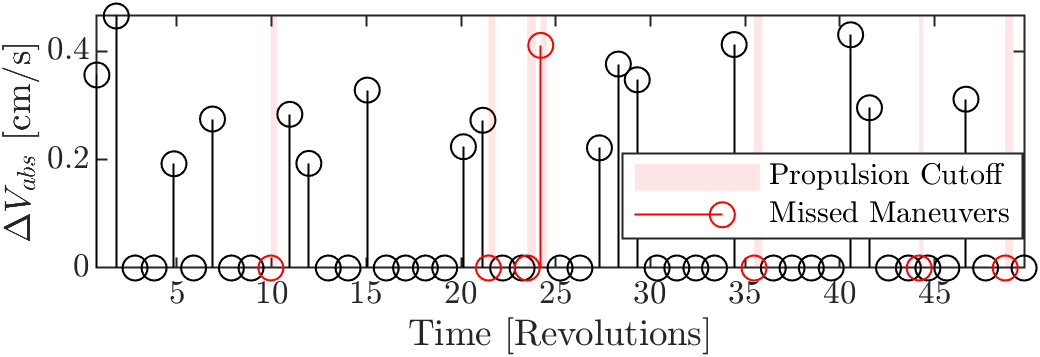}
    \end{subfigure}
    \caption{Absolute and Total $\Delta V$ Expenditure}
  \end{subfigure}
  
  \caption{Test Case 2 Trial: Five agent station-keeping for 50 revolutions}
  \label{fig:FiveAgents}
  \vspace{-.7cm}
\end{figure}

\begin{figure}[t]
  \centering
  
  \begin{subfigure}{0.59\textwidth}
    \centering
    \begin{subfigure}{0.6\textwidth}
        \centering
        \begin{subfigure}{1\linewidth}
        \centering
        \includegraphics[width=\linewidth]{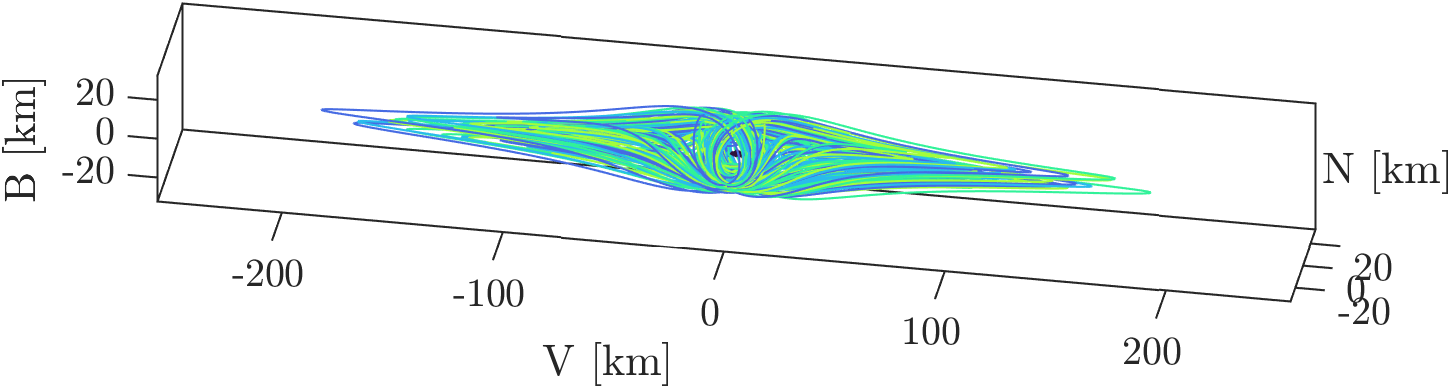}
        \vspace{.05cm}
        \end{subfigure}
        \begin{subfigure}{\linewidth}
        \centering
        \includegraphics[width=\linewidth]{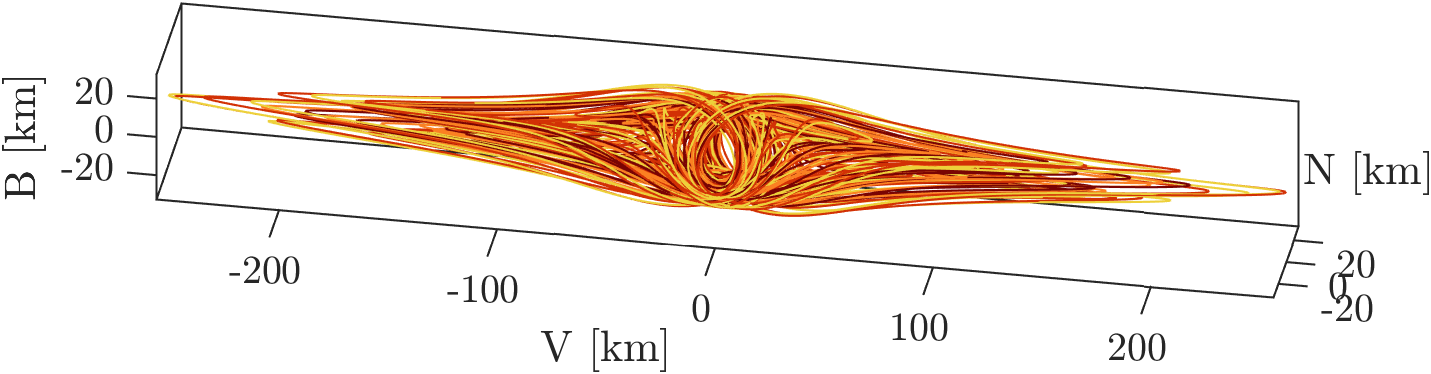}
        \vspace{.05cm}
        \end{subfigure}
        \begin{subfigure}{\linewidth}
        \centering
        \includegraphics[width=\linewidth]{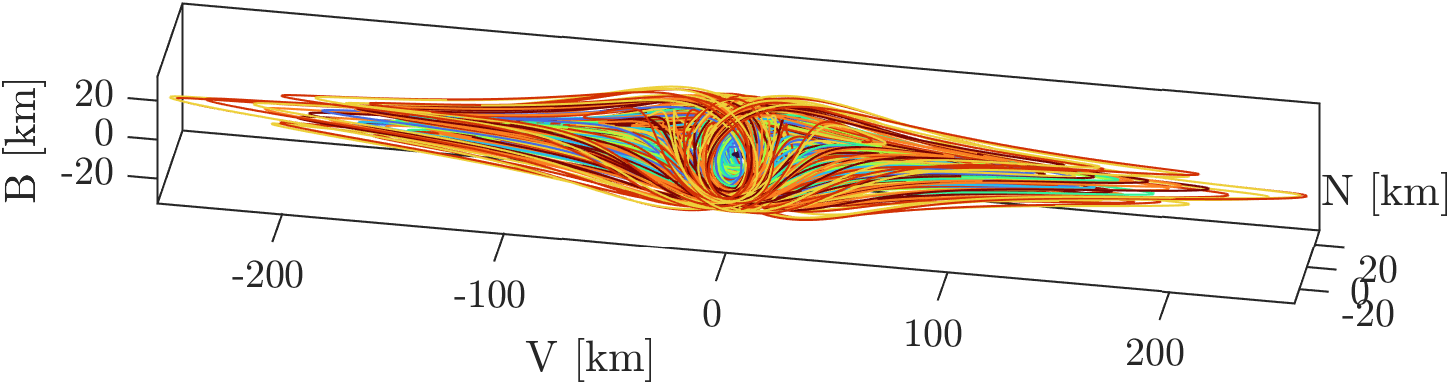}
        \end{subfigure}
    \end{subfigure}
    \hfill
    \begin{subfigure}{.39\linewidth}
      \centering
      \includegraphics[width=\linewidth]{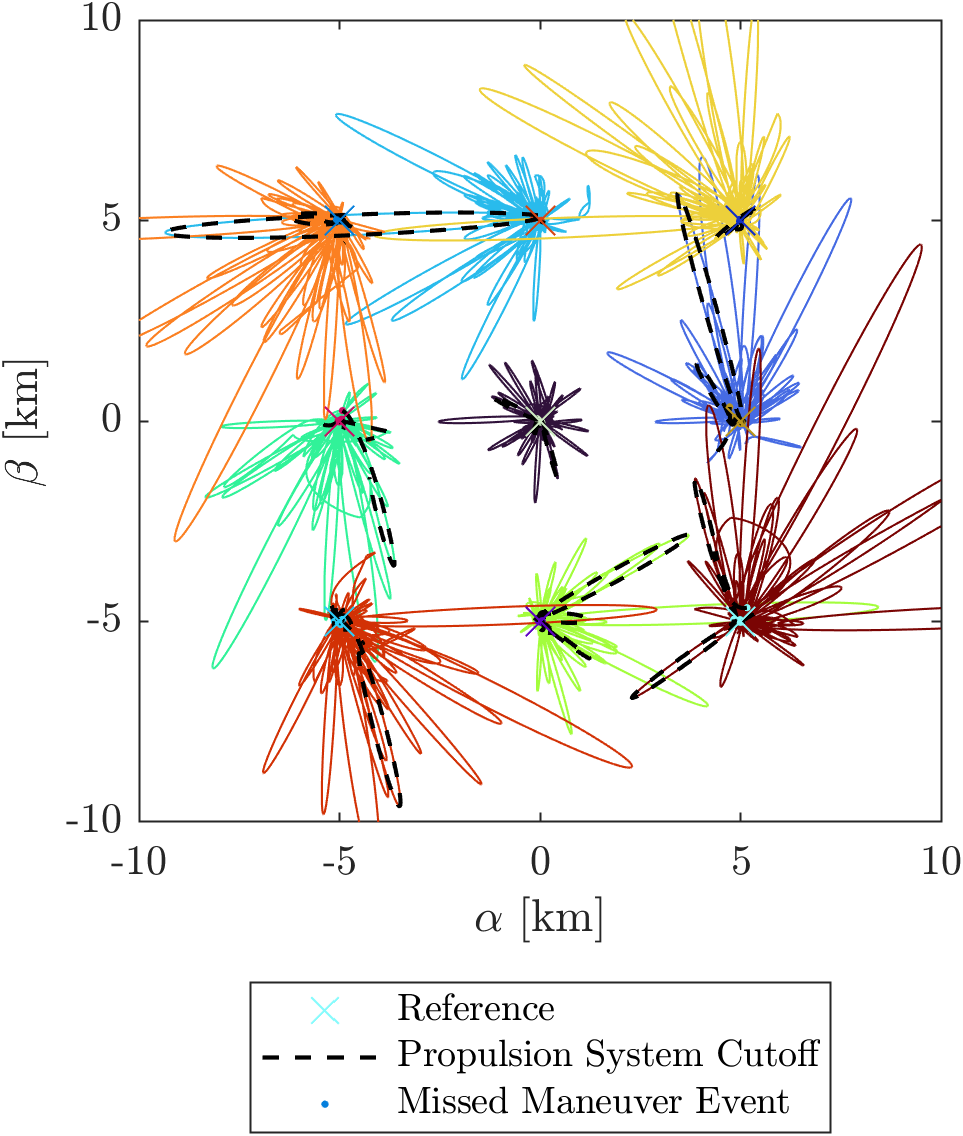}
      \vspace{.05cm}
    \end{subfigure}
    \caption{Relative Position Layers and Local Toroidal Coordinates}
  \end{subfigure}
  \hfill
  \begin{subfigure}{0.4\textwidth}
    \centering
    \includegraphics[width=\linewidth]{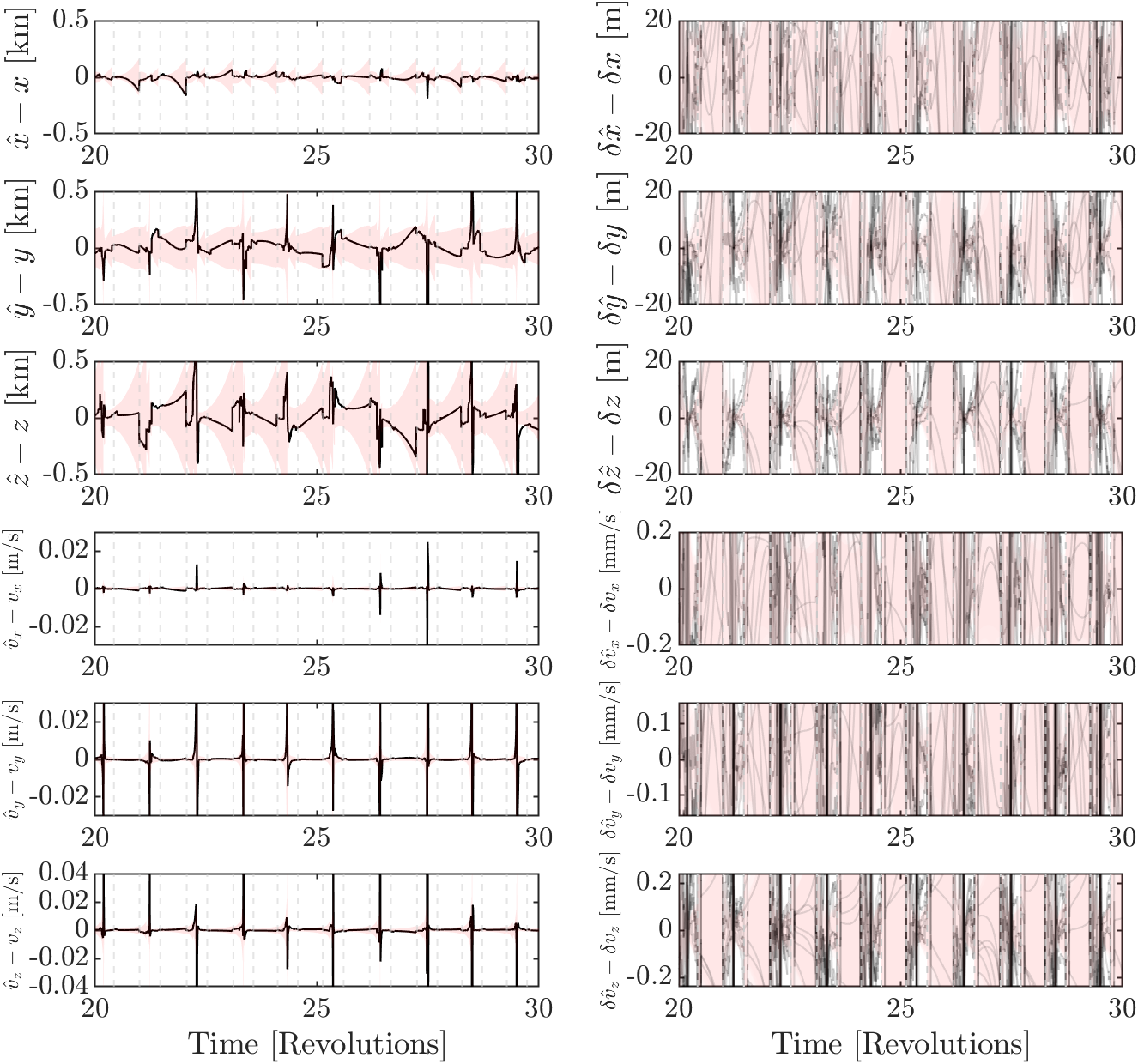}
    \caption{Navigation Subset}
  \end{subfigure}
  
  \vspace{0em}

  \begin{subfigure}{0.3\textwidth}
    \centering
    \includegraphics[width=\linewidth]{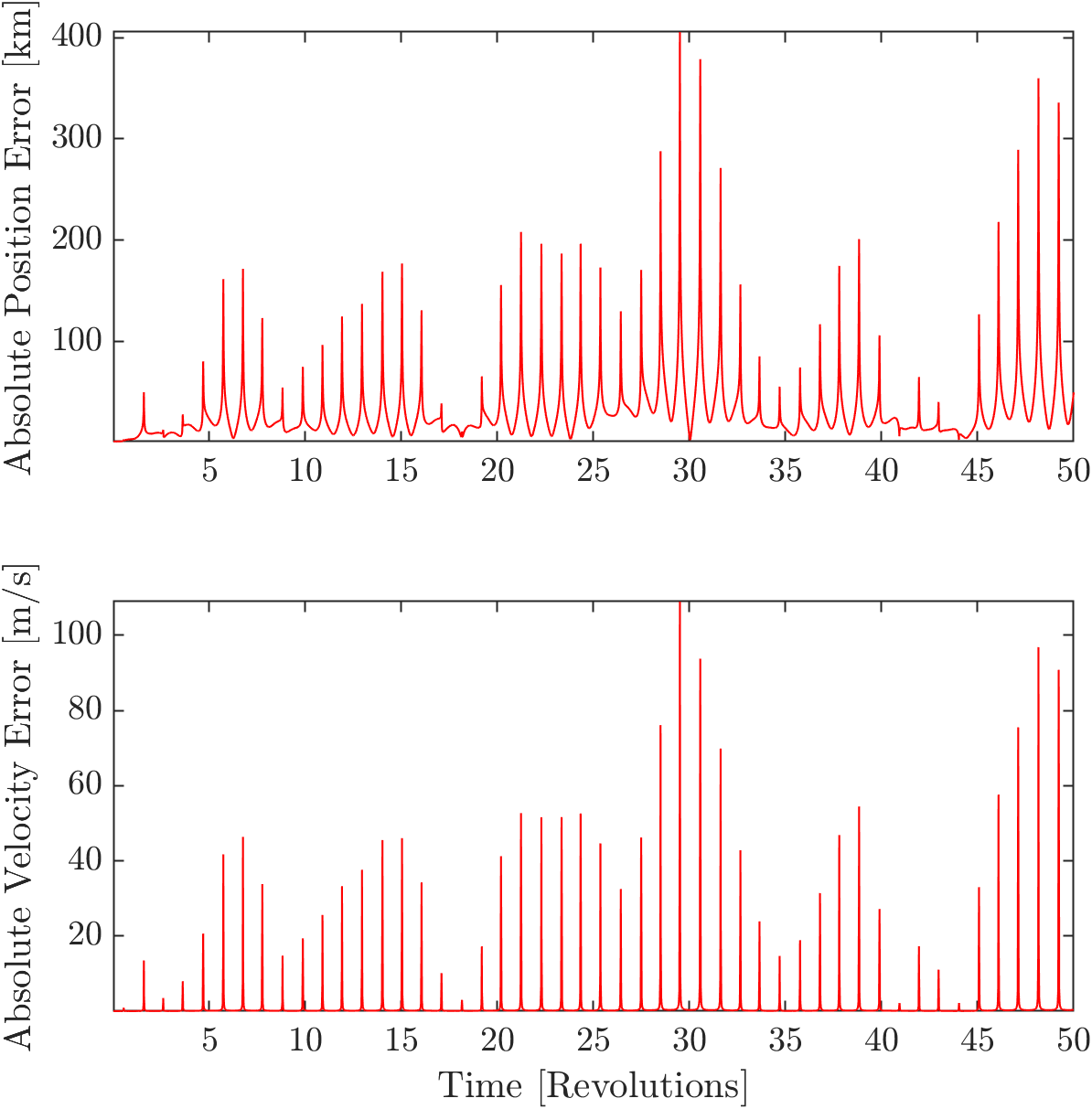}
    \caption{Absolute Reference Tracking Performance}
  \end{subfigure}
  \hfill
  \begin{subfigure}{0.3\textwidth}
    \centering
    \includegraphics[width=\linewidth]{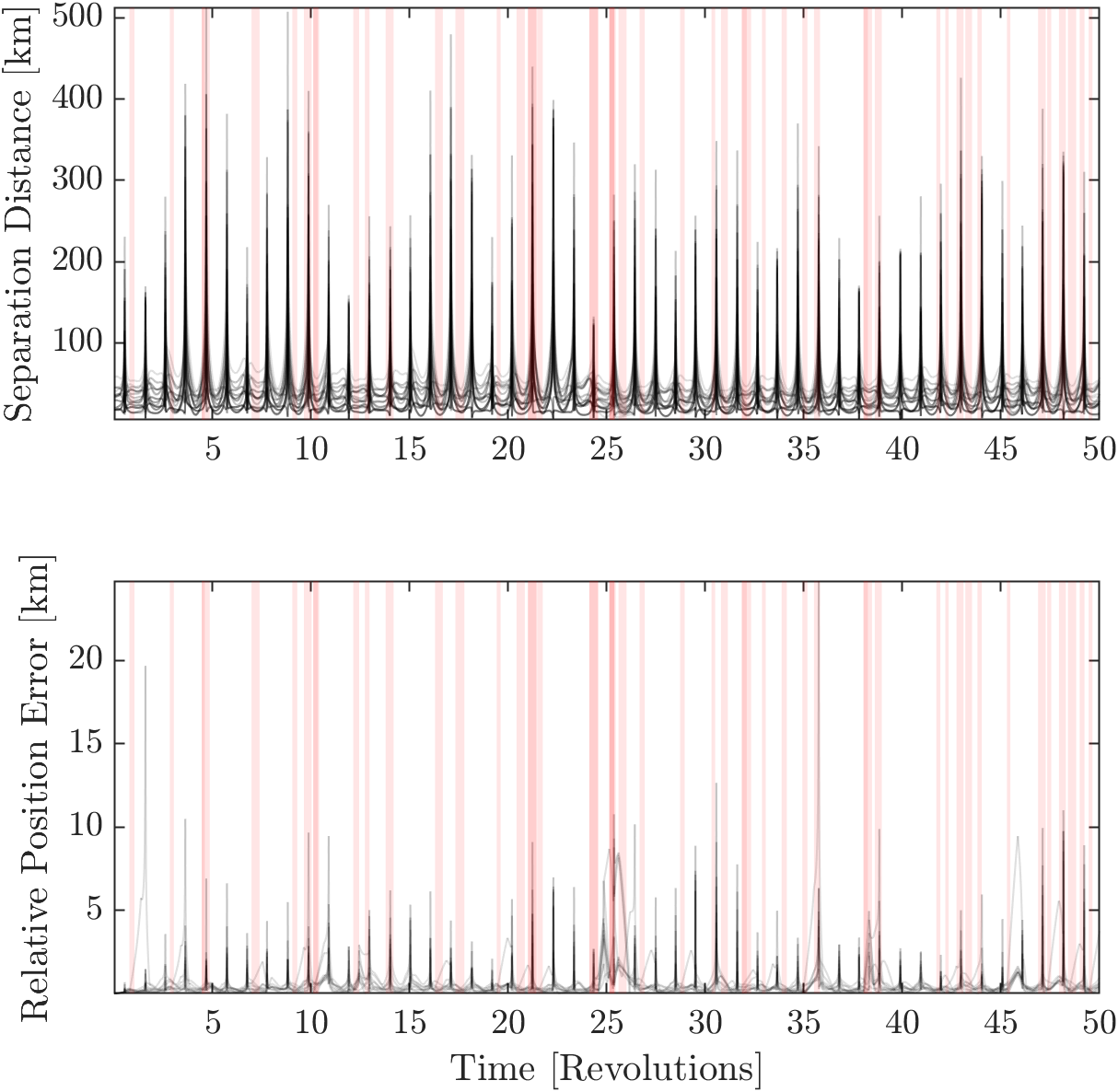}
    \caption{Relative Reference Tracking and Separation}
  \end{subfigure}
  \hfill
  \begin{subfigure}{0.32\textwidth}
    \centering
    \begin{subfigure}{.96\linewidth}
      \centering
      \includegraphics[width=\linewidth]{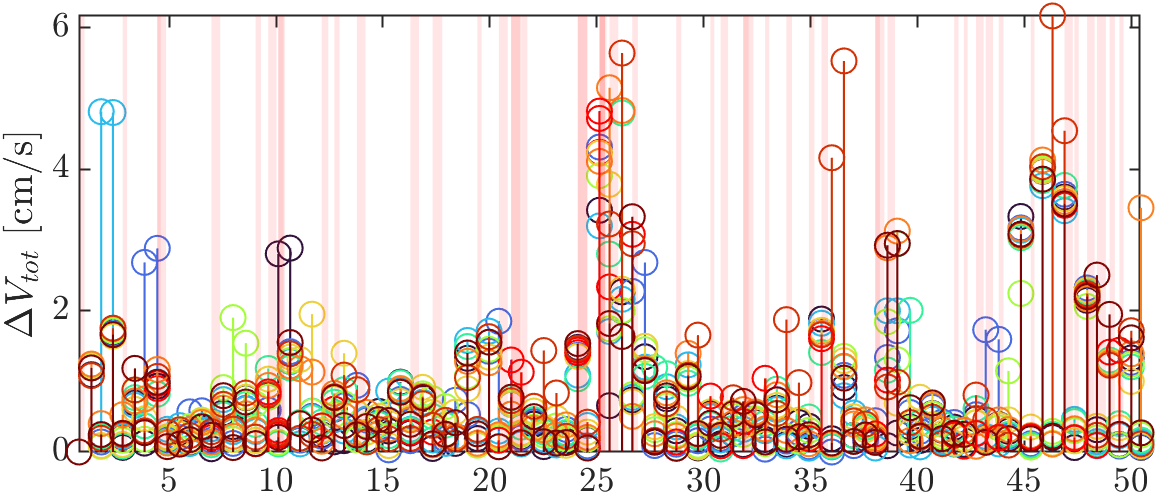}
    \end{subfigure}
    \\[1ex]
    \begin{subfigure}{\linewidth}
      \centering
      \includegraphics[width=.96\linewidth]{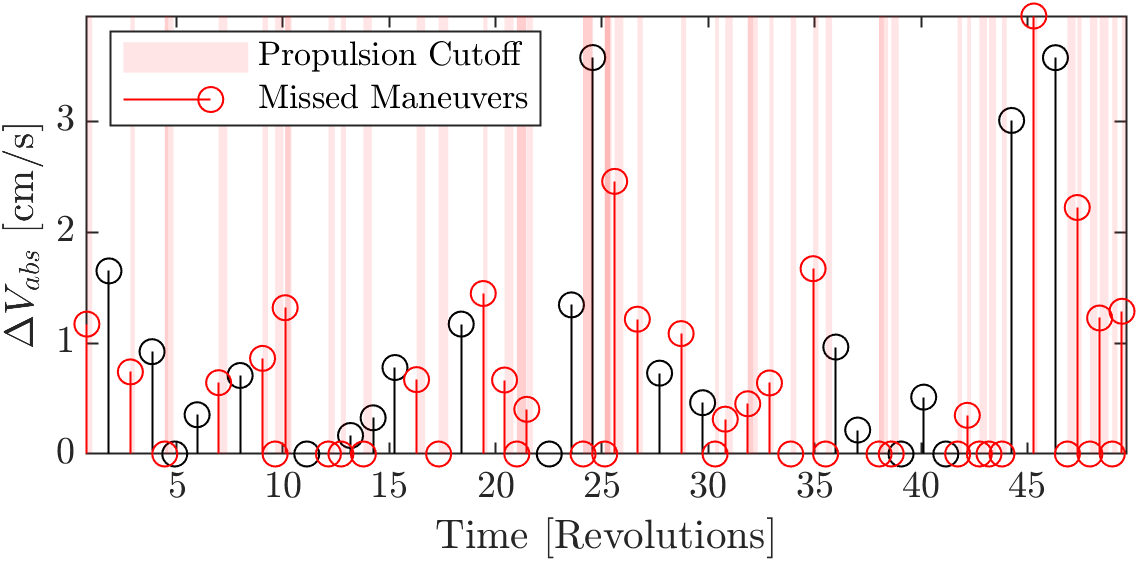}
    \end{subfigure}
    \caption{Absolute and Total $\Delta V$ Expenditure}
  \end{subfigure}
  
  \caption{Test Case 3 Trial: Nine agent station-keeping for 50 revolutions}
  \label{fig:NineAgents}
  \vspace{-.7cm}
\end{figure}

\noindent Figs.~\ref{fig:TwoAgents}, \ref{fig:FiveAgents}, and \ref{fig:NineAgents} show the performance of the proposed strategy for each respective test case, including navigation, relative and absolute tracking, $\Delta V$ histories, and trajectories in the VNB frame and local toroidal coordinates. Several performance metrics are summarized in Tab.~\ref{tab:metrics}, including yearly $\Delta V$ expenditure, separation ranges, and computation times\footnote{Experiments were performed on an Intel i7‑10750H CPU. The code was implemented and executed in MATLAB R2024a. The algorithms made use of MOSEK \cite{mosek} to optimize and YALMIP \cite{Lofberg2004} to model.}\footnote{Similar absolute control strategy implementations (e.g., Ref.~\citenum{shimane2025output}) have achieved lower computation times.}. Relative separation ranges denote the absolute minimum and maximum separation between the agents in the formation over all trials. The activation rate of the CBF-SIP represents the percentage of times that at least one iteration of the CBF optimization routine occurs, i.e. the desired $\Delta V$ is not already passively-safe.

\section{Discussion}

For all trials, the strategy demonstrates successful station-keeping for the formations with low to moderate $\Delta V$ expenditure, adequate reference tracking, low computational cost, and guarantees of passive safety. Over the cumulative 25 years of simulated absolute and relative station-keeping, the absolute maneuver and CBF-SIP schemes never fail to converge and absolute and relative bounded motion is always achieved. Moreover, no significant collision risk is detected, with the absolute lowest separation being 99 meters for Test Case 2 and 1040 meters for Test Case 1.

\subsection{Absolute Control Performance}

Absolute reference tracking is achieved on the level of hundreds of kilometers and absolute station-keeping costs are on the order of tens of centimeters per second. These figures match the findings of other studies of absolute NRHO control \cite{davis2022orbit,shimane2025output}. Lower computed absolute station-keeping maneuver magnitudes could be attributable to the connected navigation architectures, which likely produce more accurate navigation estimates. This demonstrates an important benefit associated with cislunar formation flight, whereby multiple sources of absolute measurements, when combined with relative measurements, produce improved absolute state knowledge. Previous studies, such as Refs.~\citenum{davis2022orbit,shimane2025output}, have found that absolute state knowledge error is the primary driver of $\Delta V$ consumption and absolute reference tracking, suggesting that distributed navigation architectures over cislunar formations can provide benefits to absolute control of NRHOs. A decrease in computed absolute maneuver expenditure could be additionally attributable to offloading to relative station-keeping, which cleans up after missed absolute maneuvers.

\subsection{Relative Control Performance}

The relative control concept, which leverages STMs generated about the reference trajectory for dynamics representation, QPROs as relative references, and a simple control law, demonstrates effective relative station-keeping under high-fidelity cislunar simulation. These features make the relative control law highly desirable in comparison to computing absolute maneuvers for each agent, which is incapable of providing precise relative positioning.

The formation experiences the largest relative dispersions at perilune, the region corresponding to the greatest relative separations and most chaotic dynamics. Even under low levels of noise, it is difficult to ensure good tracking at perilune as a result of only performing maneuvers at apolune. This can be compared with the operations of the PROBA-3 spacecraft \cite{Garcia2017_PROBA3_FDS}, which operates in a highly eccentric orbit about the Earth and performs high accuracy relative control only at apogee. The greatest contributor to this relative dispersion under normal conditions is the dynamical infeasibility of the QPRO. Since the QPRO is generated about a reference trajectory, and absolute control induces hundreds of kilometers of dispersion from this reference, the linearized dynamics from which the QPRO is generated degrade in accuracy. This fact has not been appropriately considered in previous studies of cislunar formation flight, which typically use simplified models and assume perfect absolute reference knowledge and tracking. Another implication of this finding is the degradation in the validity of the LTCs over the noisy, sensitive, and dispersed environment of the 9:2 NRHO.

Another source of relative tracking error is the existence of MMEs. A singular MME can severely perturb the formation geometry and is therefore the greatest source of relative disturbance. The effect of this disturbance is dependent on the size of the underlying formation, with MMEs being especially pernicious to tight formations, as evidenced by Test Case 2. The relative dispersions introduced by MMEs also contribute to increased $\Delta V$ expenditure, with the largest MMEs incurring the greatest relative $\Delta V$ to clean up.

The choice of virtual chief performs well at distributing the $\Delta V$ expenditure among the agents in both the first and second test cases. In the third test case, $\Delta V$ expenditure is less consistent across the formation. This can be attributed to the higher navigation error of the spoke spacecraft under the hub-spoke measurement architecture and the higher dispersion from the absolute reference.

Relative maneuver expenditure is also strongly correlated with the size of the formation, as seen with Test Case 3, where hundreds of kilometers of separation incur the greatest $\Delta V$ expenditure. This can be attributed to the fact that the linearized dynamics, which are used by the relative control law, become less accurate at higher separations. 



\subsection{Safety Filter Performance}

\begin{figure}[t]
    \centering
    \begin{subfigure}{\linewidth}
        \centering
        \begin{subfigure}{.5\linewidth}
            \centering
            \includegraphics[width=\linewidth]{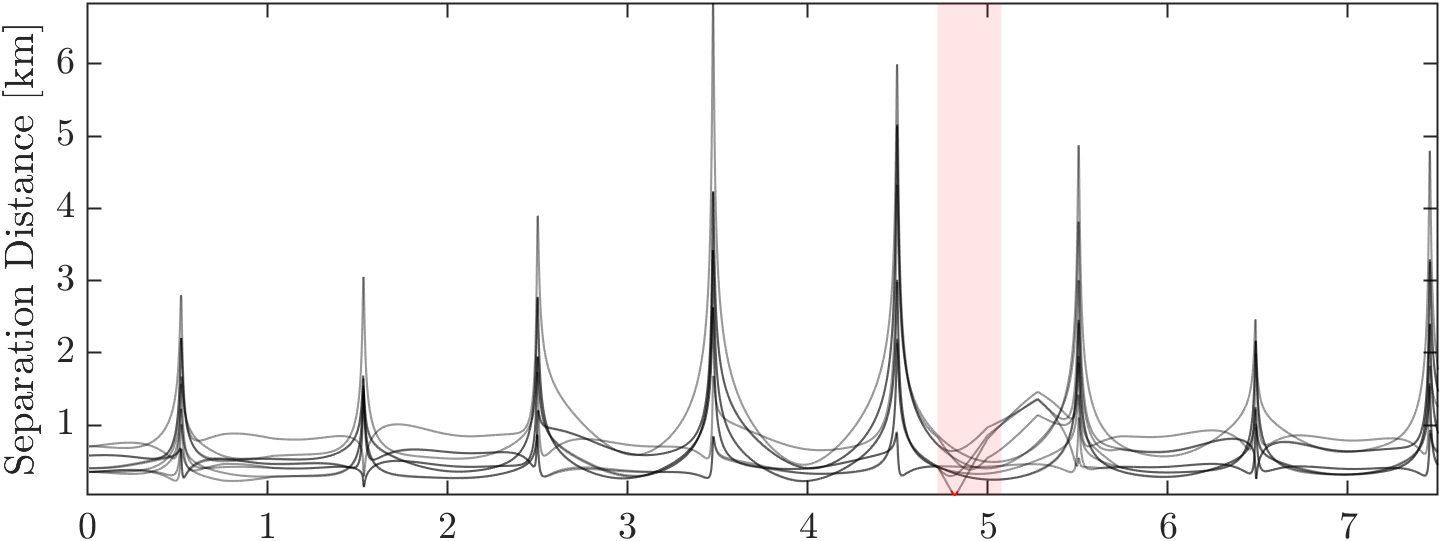}
        \end{subfigure}%
        \hspace{0.5em}
        \begin{subfigure}{.34\linewidth}
            \centering
            \includegraphics[width=\linewidth]{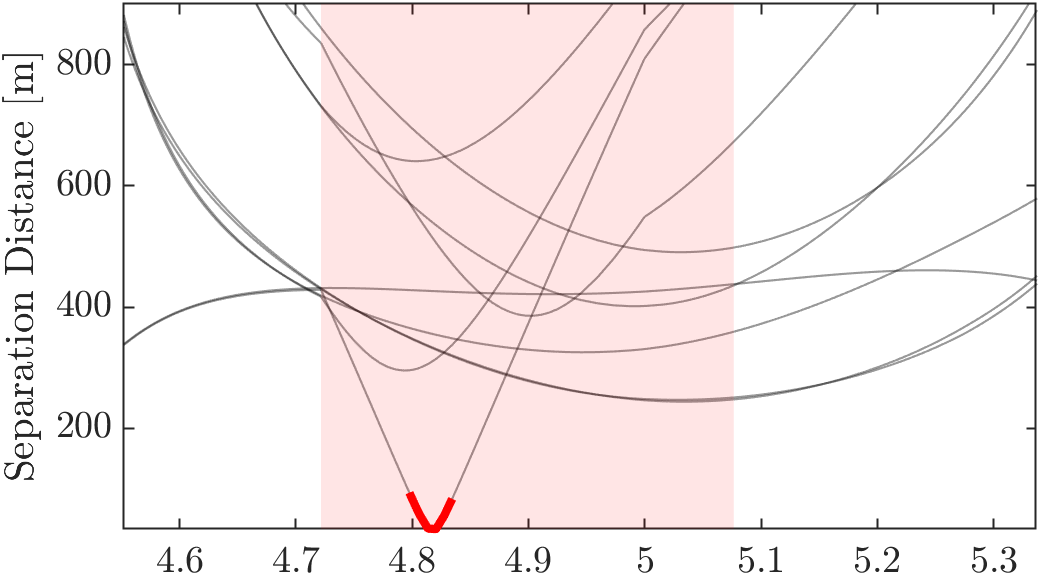}
        \end{subfigure}
        \caption{CBF-SIP not applied}
        \vspace{-.2cm}
    \end{subfigure}

    \vspace{0.5em} 

    \begin{subfigure}{\linewidth}
        \centering
        \begin{subfigure}{.5\linewidth}
            \centering
            \includegraphics[width=\linewidth]{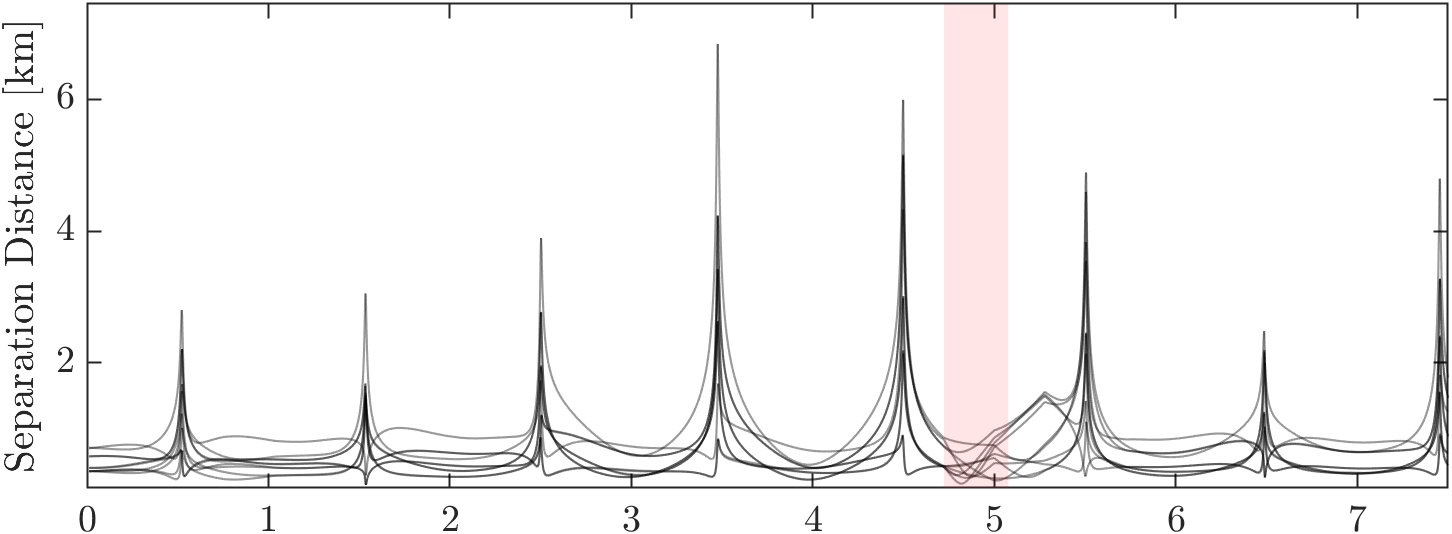}
        \end{subfigure}%
        \hspace{0.5em}
        \begin{subfigure}{.34\linewidth}
            \centering
            \includegraphics[width=\linewidth]{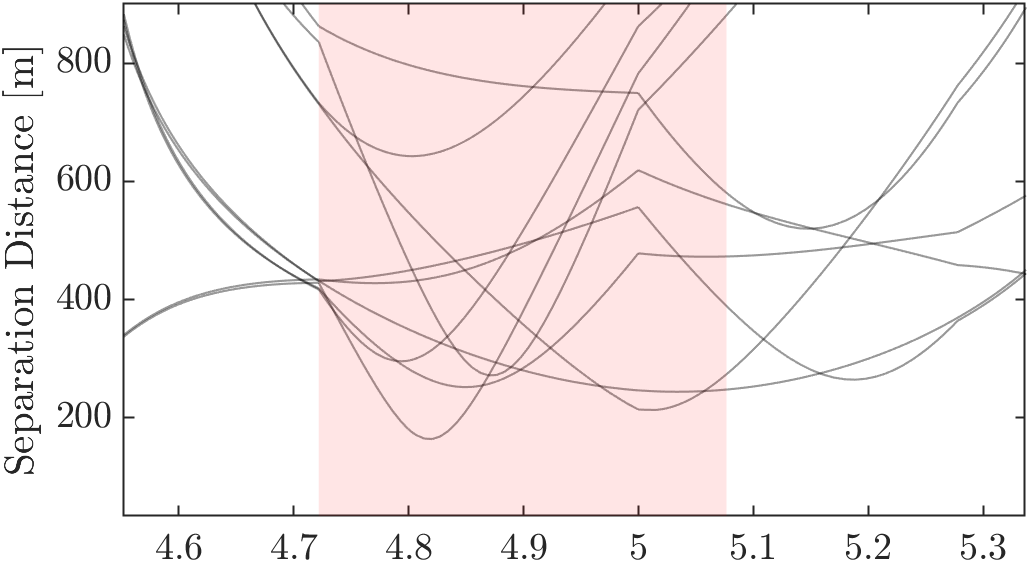}
        \end{subfigure}
        \caption{Active CBF-SIP}
        \vspace{-.2cm}
    \end{subfigure}

    \caption{Five agent near miss scenario}
    \label{fig:safety}
    \vspace{-.7cm}
\end{figure}

A key feature of the proposed CBF-SIP is that optimization routines are only called if the formation enters an unsafe condition, leveraging the rarity of near-miss collision events. As a result, Test Case 3 features almost no activations of the CBF-SIP, meaning runtime is only spent verifying that ensuing trajectories are safe. In contrast, Test Case 2 features a 29\% activation rate of the CBF-SIP, given the tighter formation, and yet retains long-term stability over the duration of flight. 

The ability of the CBF-SIP to guarantee safety is demonstrated in Fig.~\ref{fig:safety}, which originates from Test Case 2. When the CBF-SIP is not active, a near-miss event occurs with the spacecraft trajectories coming to within 38 meters of each other following a missed station-keeping maneuver. With the CBF-SIP active, the CBF appropriately corrects the maneuvers such that the near-miss event is avoided and the closest approach occurs at 100 meters, the prescribed bound. 

\section{Conclusion}

This work presents a guidance and control strategy for long duration station-keeping in cislunar near-rectilinear halo orbits and demonstrates safe and efficient deployment in high-fidelity scenarios. The strategy represents one of the first demonstrations of a complete GN\&C architecture applied to cislunar formation flight in NRHOs under high-fidelity simulation. The combination of state-of-the-art absolute reference tracking strategies with a simple control law for QPRO tracking is achieved through the introduction of a control barrier function safety filter, which allows for decoupling of the absolute and relative control formulations. Relative navigation, relative control, and the CBF-SIP, are computationally efficient and introduce minimal additional computational burden, meaning that these aspects of formation-keeping could be performed autonomously. In this respect, the proposed framework has the potential to introduce minor amounts of additional operational burden compared to a single spacecraft mission.

Finally, it is emphasized that though this work has focused on the relevant application of the 9:2 $L_2$ NRHO, the proposed strategy is generalizable to other multi-body orbits. Future work includes the systematic verification of the GN\&C strategy to other orbits of interest, including distant retrograde orbits, butterfly orbits, Lyapunov orbits, and other $L_1/L_2$ Halo orbits, the inclusion of advanced reconfiguration strategies, such as those proposed in Refs.~\citenum{takubo2025safeoptimalnspacecraftswarm,khoury_relative_2024}, and further analysis of strategies for generating longer-duration QPROs. In addition, the inclusion of greater sources of noise, such as larger solar radiation pressure effects and attitude desaturation events, would be beneficial to include for further validation.

\section{Acknowledgements}

This work is supported by Blue Origin (SPO \#299266) as an Associate Member and Co-Founder of Stanford’s Center of AEroSpace Autonomy Research (CAESAR). This article solely reflects the opinions and conclusions of its authors and not any sponsors.

\bibliographystyle{AAS_publication}   
\bibliography{references}   

\end{document}